\documentclass[11pt,twoside,reqno]{amsart}

\usepackage{bm}
\usepackage{amssymb,amsmath,amstext,amsthm,amsfonts,amscd,xcolor}
\usepackage[a4paper, left=3cm, right=3cm, top=3cm, bottom=3cm]{geometry}
\usepackage{dsfont}
\usepackage{amsthm, enumerate}
\usepackage[ansinew]{inputenc}
\usepackage{graphicx}
\usepackage[mathscr]{eucal}
\usepackage{xcolor} 
\usepackage{mathrsfs}
\usepackage{hyperref}
\makeatletter
\newsavebox\myboxA
\newsavebox\myboxB
\newlength\mylenA

\usepackage{fancyhdr}
\numberwithin{equation}{section}
\newtheorem{theorem}{Theorem}[section]
\theoremstyle{definition}
\newtheorem{definition}{Definition}[section]  
\newtheorem{lemma}{Lemma}[section]      
\newtheorem{corollary}{Corollary}[section]    
\newtheorem{proposition}{Proposition}[section] 
\theoremstyle{remark}
\newtheorem{remark}{Remark}[section]
\title{Abstract Continuity Theorem for the Lyapunov Exponents of linear cocycles}

\begin{document}
	
	\author[A. Cai]{Ao Cai}
	\address{School of Mathematical Sciences, Soochow University, Soochow 215006, China}
	\email{acai@suda.edu.cn}
	
		\author[X. Deng]{Xiaojuan Deng}
	\address{School of Mathematical Sciences, Soochow University, Soochow 215006, China}
	\email{xiaojuanky@163.com}

\begin{abstract}
	We prove H\"older continuity of Lyapunov exponents for general linear cocycles varying measures on the base in Wasserstein distance under the assumption of uniform large deviations type (LDT) estimates. This is a measure version of the abstract continuity theorem (ACT) established by Duarte-Klein \cite{duarte2016lyapunov}. The main obstacle here lies in the fact that the magnitude of the exceptional sets in LDT estimates is constantly changing when the base measures deviate. We overcome this via a combination of a Urysohn-type lemma and properties of Wasserstein distance in every iteration step. Our measure version of ACT, combined with the original work of Duarte-Klein \cite{duarte2016lyapunov}, provides a complete scheme for proving joint H\"older continuity of Lyapunov exponents with respect to both measure and fiber which resolves all parameter dependence. This continuity theorem is general and applicable to a wide range of mathematical models, including product of random matrices and cocycles essentially generated by shifts. In particular, it applies to associated Schr\"odinger operators which are  central objects in the study of mathematical physics.
\end{abstract}

\maketitle

\section{Introduction}\label{sec:intro}
An ergodic measure preserving dynamical system, MPDS for short, $(X,\mathcal{F},\mu,T)$ consists of a Polish metric space $X$, a $\sigma$-algebra $\mathcal{F}$, a probability measure $\mu$ defined on $(X,\mathcal{F})$ and a transformation $T:X\to X$ that is both ergodic and measure preserving.
\par Given a measurable function $A:X\to \mathrm{Mat}(d,\mathbb{R})$, a linear cocycle over an ergodic MPDS $(X,\mathcal{F},\mu,T)$ is a skew-product map $F$ defined by
$$F: X\times\mathbb{R}^d\ni (x,v)\mapsto(Tx,A(x)v)\in X\times\mathbb{R}^d.$$
We identify $F$ with the pair $(\mu,A)$ and call $T$ the base dynamics while $A$ defines the fiber action. The iterates of the cocycle $F$ are given by $(T^nx,A^{(n)}(x)v)$, where
$$A^{(n)}(x)=A(T^{n-1}x)\cdot\cdots\cdot A(Tx)\cdot A(x).$$
We say $(\mu,A)$ is a continuous cocycle if it satisfies that every $A:X\to \mathrm{Mat}(d,\mathbb{R})$ is a continuous function where $d\in\mathbb{N}$ and $\mathrm{supp}(\mu)$ is contained in a compact set. We denote the space of continuous cocycles by $\mathscr{C}$. Also, let $P(X)$ denote the space of Borel probability measures on $X$.
\par The Lyapunov exponents of the cocycle $(\mu,A)$ are numbers $L_1(\mu,A)\geq L_2(\mu,A)\geq\cdots\geq L_d(\mu,A)\geq-\infty$, which measure the average exponential growth of iterates of the cocycle along invariant fiber subspaces. By Kingman's sub-additive ergodic theorem, these Lyapunov exponents can be characterized by the $\mu$-almost everywhere limits for $1\leq k\leq d$,$$L_k(\mu,A)=\lim_{n\to +\infty}\frac{1}{n}\log s_k(A^{(n)}(x)),$$
where $s_1(A^{(n)}(x))\geq s_2(A^{(n)}(x))\geq\cdots\geq s_d(A^{(n)}(x))\geq 0$ denote the singular values of the matrix $A^{(n)}(x)\in \mathrm{Mat}(d,\mathbb{R})$.
\par In this paper, we are interested in the continuity properties of the Lyapunov exponents of general linear cocycles with respect to measure $\mu$. By fixing the fiber action and considering different measures on the base that are sufficiently close in the Wasserstein distance $W_1$ elaborated in the preliminaries, we demonstrate the continuity of the Lyapunov exponents and show that the Lyapunov blocks are H\"older continuous when the uniform LDT estimates (in measure) and the gaps between Lyapunov exponents are available. In this sense, our ACT bridges the LDT in Probability Theory to the continuity of Lyapunov exponents in Dynamical System. As the continuity of Lyapunov exponents is a classic and challenging problem in history, our ACT reduces it to proving uniform LDT estimates which are hard yet operable for various systems in practice.
\par Large deviation principle primarily investigates rare events and their asymptotic behavior on an exponential scale. The fundamental idea of large deviation theory was first proposed by H. Cram\'er in 1938 \cite{cramer1994collected}: for a sequence of independent and identically distributed random variables, the tail probability of their partial sums decays at an exponential rate. The exact statement of the theorem can be found in Dembo \cite{dembo2009large}: let $X_0,X_1,X_2,\cdots$ be a real-valued random process and denote by $S_n=\sum^{n-1}_{j=0}X_j$ the corresponding sum process. If the random process $\{X_n\}$ is i.i.d. with mean $\mathbb{E}(X_0)$ and finite moment generating function $M(t):=\mathbb{E}[e^{tX_0}]<+\infty$ for all $t>0$, then $$\lim_{n\to+\infty}\frac{1}{n}\log\mu_n[|\frac{1}{n}S_n-\mathbb{E}(X_0)|>\varepsilon]=-I(\varepsilon),$$ where $I(\varepsilon):=\sup_{t>0}(t\varepsilon-\log M(t)+t\mathbb{E}(X_0))$ is called the rate function.
\par The precise formulation of the large deviation principle is introduced by Varadhan \cite{varadhan1966asymptotic}: for a lower semicontinuous function $I:\mathbb{R}\to [0,\infty]$, a positive increasing sequence $a_n\to\infty$, and a sequence of measures $\{\mu_n\}\subset P(X)$, we have:
\begin{align}
    & \limsup_{n\to\infty}\frac{1}{a_n}\log\mu_n(F)\leq-\inf_{x\in F}I(x)\qquad\forall \text{ closed }  F\subset X,\\
    &\liminf_{n\to\infty}\frac{1}{a_n}\log\mu_n(G)\geq-\inf_{x\in G}I(x)\qquad\ \forall \text{ open } G\subset X.
\end{align}
Later, Donsker and Varadhan \cite{donsker1975asymptotic1,donsker1975asymptotic2,donsker1976asymptotic3,donsker1983asymptotic4} showed the renowned Donsker-Varadhan theoretical framework, which quantifies the probability of large deviations events through the definition of a rate function. For detailed discussions of other large deviation principles, including those for Markov processes, we refer to \cite{rassoul2015course}.
\par Given a linear cocycle $A:X\to \mathrm{Mat}(d,\mathbb{R})$ over the ergodic MPDS $(X,\mathcal{F},\mu,T)$, we say that $\mu$ satisfies a large deviation principle with rate function $I(\varepsilon)$ if
$$\lim_{n\to+\infty}\frac{1}{n}\log\mu\{|\frac{1}{n}\log\|A^{(n)}(x)\|-L_1(\mu)|>\varepsilon\}=-I(\varepsilon).$$
\par However, our study of the continuity of Lyapunov exponents for linear cocycles only requires good upper bounds on the measure of tail events. 
\par In our work, we require the precision of large deviations type estimates and the exponential decay of the deviation measure (sub-exponential decay is allowed but the result shall be reasonably weaker, so we do not elaborate on this). We use the rate function $r(n) = e^{-cn}$ to describe the measure of the set of $x$, for which the deviation of the average of the sum of random variables from the mean is at least $\varepsilon$. Here, we refer to this set as the deviation set and the constant $c>0$ is related to $\varepsilon > 0$. Let $\mathcal{J}$ be the space of all rate functions $r(n)$, we specify that $\mathcal{J}$ is convex, i.e. $r(n)$ is closed under scalar multiplication and addition within $\mathcal{J}$. 
\par Next, consider $\mathcal{C}^\alpha(X)$ as the space of all $\alpha$-H\"older continuous functions $\xi:X\to\mathbb{R}$, where $0<\alpha\leq 1$. We define some good upper bounds on the measure of tail events as Large Deviations Type (LDT) estimates in the following.
\begin{definition}\label{de1}
    An observable $\xi\in \mathcal{C}^\alpha(X)$ satisfies a base-LDT estimate if for all $\varepsilon>0$, there are $\bar{m}=\bar{m}(\xi,\varepsilon)\in\mathbb{N},c=c(\xi,\varepsilon)>0$ such that for all $n\geq\bar{m}$, we have 
    \begin{equation}
        \mu\{x\in X:|\frac{1}{n}\sum_{j=0}^{n-1}\xi(T^jx)-\int_X\xi d\mu(x)|>\varepsilon\}<r(n)= e^{-cn}.
    \end{equation}
\end{definition}
\begin{definition}\label{de2}
    An observable $\xi\in \mathcal{C}^\alpha(X)$ satisfies a uniform base-LDT estimate if for all $\varepsilon>0$, there are $\bar{\delta}=\bar{\delta}(\xi,\varepsilon,\mu)>0,\bar{m}=\bar{m}(\xi,\varepsilon,\mu)\in\mathbb{N}$ and $c'=c'(\xi,\varepsilon,\mu)>0$ such that for all $\nu\in P(X)$ with $W_1(\nu,\mu)<\bar{\delta}$ and for all $n\geq\bar{m}$, we have 
    \begin{equation}
        \nu\{x\in X:|\frac{1}{n}\sum_{j=0}^{n-1}\xi(T^jx)-\int_X\xi d\nu(x)|>\varepsilon\}<r'(n)= e^{-c'n}.
    \end{equation}
\end{definition}
\begin{definition}\label{de3}
    A continuous cocycle $(\mu,A)\in\mathscr{C}$ satisfies a fiber-LDT estimate if for all $\varepsilon>0$, there are $\bar{m}=\bar{m}(A,\varepsilon)\in\mathbb{N},c=c(A,\varepsilon)>0$ such that for all $n\geq\bar{m}$, we have 
    \begin{equation}
        \mu\{x\in X:|\frac{1}{n}\log\|A^{(n)}(x)\|-L_1^{(n)}(\mu)|>\varepsilon\}<r(n)= e^{-cn}.
    \end{equation}
\end{definition}
\begin{definition}\label{de4}
    A continuous cocycle $(\mu,A)\in\mathscr{C}$ satisfies a uniform fiber-LDT estimate if for all $\varepsilon>0$, there are $\bar{\delta}=\bar{\delta}(A,\mu,\varepsilon)>0,\bar{m}=\bar{m}(A,\mu,\varepsilon)\in\mathbb{N}$ and $c''=c''(A,\mu,\varepsilon)>0$ such that for all $\nu\in P(X)$ with $W_1(\nu,\mu)<\bar{\delta}$ and for all $n\geq\bar{m}$, we have 
    \begin{equation}
        \nu\{x\in X:|\frac{1}{n}\log\|A^{(n)}(x)\|-L_1^{(n)}(\nu)|>\varepsilon\}<r''(n)= e^{-c''n}.
    \end{equation}
\end{definition}
\par Under the premise that the uniform LDT estimates hold, we can prove our main theorem as follows.
\begin{theorem}\label{thm1}
   Given an ergodic MPDS $(X,\mu,\mathcal{F},T)$, a space of continuous cocycle $\mathscr{C}$, we assume that:\\
   (1) every $\xi\in \mathcal{C}^\alpha(X)(0<\alpha\leq 1)$ satisfies a uniform base-LDT estimate.\\
   (2) every cocycle $(\mu,A)\in\mathscr{C}$ for which $L_1(\mu)>L_2(\mu)$ satisfies a uniform fiber-LDT estimate.
   \par Then there exists a neighborhood $\mathcal{V}_\mu$ of $\mu$, such that for any ergodic probability measure $\nu\in\mathcal{V}_\mu$, the first Lyapunov exponent $L_1(\nu)$ is H\"older continuous with respect to $\nu$ in Wasserstein distance. The modulus of continuity is given by $\omega(h):=[r(a\log \frac{1}{h})]^{\frac{1}{2}-}$ for some $r=r(\mu)\in\mathcal{J}$ and $a = a(\mu)> 0$, where $\frac{1}{2}-$ represents a value slightly less than $\frac{1}{2}$.
\end{theorem}  

\begin{remark}
	The assumption of $A$ being continuous is necessary due to comparison of the exceptional sets under distinct measures in each iteration step. This is essentially different from the fiber version of ACT where the base measure is always fixed. Moreover, the requirement of base-LDT estimates to be uniform is inevitable as base measures vary, which is also diverse from Duarte-Klein \cite{duarte2016lyapunov}.
\end{remark}
\begin{remark}\label{re2}
    If we assume $L_k(\mu)>L_{k+1}(\mu)$, then locally near $\mu$ the Lyapunov blocks $L_1+L_2+\cdots+L_k$ are H\"older continuous with respect to $\nu$ where the modulus of continuity $\omega(h)$ is given in Theorem \ref{thm1}. In this case, the continuity of each $L_i(1\leq i\leq k)$ can also be obtained, but when $k<i\leq d$, we cannot determine the continuity of $L_i$. However, if we assume $A$ is invertible, we can prove all Lyapunov exponents are continuous. The requirement for invertibility of $A$ seems necessary to prove the continuity of all Lyapunov exponents in the Wasserstein distance. Otherwise, even slight perturbations of the measure could cause the Lyapunov exponents to drop directly from a finite number to negative infinity essentially because a set of $\mu$-measure zero may no longer retain zero in $\nu$-measure under perturbation.
\end{remark}
In connection with Theorem 3.1 in Duarte-Klein \cite{duarte2016lyapunov}, and assuming that uniform (in measure and fiber) LDT estimates are satisfied simultaneously, we can obtain joint H\"older continuity of Lyapunov exponents in measure and fiber. \par
We assume that the space of continuous matrix-valued functions defined on Polish space $X$ is endowed with the distance $$d_\infty(A,B)=\sup_{x\in X}\|A(x)-B(x)\|,$$ where any  $A(x),B(x)\in\mathrm{Mat}(d,\mathbb{R})$ are continuous. Therefore, on the space $\mathscr{C}$, we consider the product metric $d$, which is defined by $$d((\mu,A),(\nu,B))=\max\{W_1(\mu,\nu),d_\infty(A,B)\},$$ for any cocycles $(\mu,A),(\nu,B)\in\mathscr{C}$. Analogous to Definition \ref{de4}, we can define the following uniform LDT estimates.
\begin{definition}\label{de5}
	A continuous cocycle $(\mu,A)\in\mathscr{C}$ satisfies a uniform fiber-LDT estimate (in measure and fiber) if for all $\varepsilon>0$, there are $\bar{\delta}=\bar{\delta}(A,\mu,\varepsilon)>0,\bar{m}=\bar{m}(A,\mu,\varepsilon)\in\mathbb{N}$ and $c''=c''(A,\mu,\varepsilon)>0$ such that for all $(\nu,B)\in\mathscr{C}$ with $d((\mu,A),(\nu,B))<\bar{\delta}$ and for all $n\geq\bar{m}$, we have 
\begin{equation}
	\nu\{x\in X:|\frac{1}{n}\log\|B^{(n)}(x)\|-L_1^{(n)}(\nu,B)|>\varepsilon\}<r''(n)= e^{-c''n}.
\end{equation}
\end{definition}
We formulate the joint H\"older continuity result as follows:
\begin{theorem}\label{thm1.2}
	Given an ergodic MPDS $(X,\mu,\mathcal{F},T)$, a space of continuous cocycle $\mathscr{C}$ such that\\
	(1) every $\xi\in \mathcal{C}^\alpha(X)(0<\alpha\leq 1)$ satisfies a uniform base-LDT estimate.\\
	(2) every cocycle $(\mu,A)\in\mathscr{C}$ for which $L_1(\mu,A)>L_2(\mu,A)$ satisfies a uniform fiber-LDT estimate (in measure and fiber).\par
	There exists a neighborhood $\mathcal{V}$ of $(\mu,A)$ such that for all $(\nu,B)\in\mathcal{V}$, the first Lyapunov exponent $L_1(\nu,B)$ is H\"older continuous on the measure and fiber. 
\end{theorem}
Additionally, we can discuss the joint H\"older continuity of Lyapunov blocks for the case when $L_k(\mu,A)>L_{k+1}(\mu,A)$, similar to Remark \ref{re2}.
\begin{remark}
	Theorem \ref{thm1} \& \ref{thm1.2} can further demonstrate that counter examples of Bochi-Ma\~n\'e \cite{BOCHI2002} cannot simultaneously satisfy uniform base and fiber LDT estimates. More precisely, our belief is that the uniform fiber-LDT estimates break down, while the uniform base-LDT estimates could still be valid. The underlying intuition is that Lyapunov exponents are always upper semicontinuous, and uniform base-LDT implies the upper semicontinuity of Lyapunov exponents. Moreover, the non-uniform hyperbolicity in Bochi-Mañé type counter examples leads to excessive variability in the growth behavior of matrix products along different orbits, making it impossible to control them with a uniform rate function, which destroys the uniform fiber-LDT estimates.
\end{remark}
As is mentioned before, our Theorem \ref{thm1} \& \ref{thm1.2} use uniform LDT estimates as input to directly receive continuity of Lyapunov exponents as output. For the latter, it experiences a long and fruitful history of development \cite{avila2010extremal,avila2014complex,bourgain2002continuity,cai2023furstenberg,furstenbergkifer1983,goldstein2001holder,jitomirskaya2009continuity,le1989regularite,malheiro2015lyapunov}. Note that Theorem \ref{thm1} \& \ref{thm1.2} naturally include cocycles essentially generated by shifts as applicable objects (in this case ergodicity is trivial when the base measure varies), which contain but are not restricted to random cocycles and mixed random-quasiperiodic cocycles once uniform LDT estimates are available (note that the uniform LDT estimates of the latter have not been established yet).
\par The study of the continuity of Lyapunov exponents for random cocycles can be traced back to the 1980s. Furstenberg and Kifer proved in \cite{furstenbergkifer1983} that if the measure supported on a subgroup of $\mathrm{GL}(d,\mathbb{R})$ is irreducible, then the first Lyapunov exponent is continuous on the measure $\mu$. Le Page \cite{le1989regularite} assumed the strong irreducibility and contraction for $\mathrm{GL}(d,\mathbb{R})$-valued cocycles and obtained the local H\"older continuity of the first Lyapunov exponent of random matrix products with respect to the random potential. A few years ago, in Malheiro and Viana \cite{malheiro2015lyapunov}, it turns out that the extremal Lyapunov exponents of $\mathrm{GL}(2,\mathbb{R})$-valued cocycles over mixing Markov shifts depend continuously on the matrix coefficients and the transition probabilities. Similarly, Backes, Brown, and Bulter \cite{backes2018continuity} obtained that Lyapunov exponents vary continuously with respect to the fiber and measure for fiber-bunched $\mathrm{GL}(2,\mathbb{R})$-valued cocycles, which answers a conjecture of Viana. Bocker-Neto and Viana \cite{bocker2017continuity} showed that the extremal Lyapunov exponents are jointly continuous regarding the two-dimensional fiber $A$ and the probability weights $p$ for locally constant $\mathrm{GL}(2,\mathbb{C})$-valued cocycles over Bernoulli shifts. Subsequently, Avila, Eskin, and Viana \cite{avila2023continuity} extended this result, proving that for random $\mathrm{GL}(d,\mathbb{C})$-valued cocycles (or $\mathrm{GL}(d,\mathbb{R}),d\geq 2$), all Lyapunov exponents are continuous with respect to joint parameters of the fiber $A$ and the probability weights $p$. Moreover, Tall and Viana \cite{tall2020moduli} proved that for random $\mathrm{GL}(2,\mathbb{R})$ cocycles the Lyapunov exponents are point-wisely Hölder continuous on measure in Wasserstein distance when the Lyapunov exponents are distinct and they are log-H\"older continuous at every point.
\par In recent years, Duarte and Klein have been dedicated to studying the properties of Lyapunov exponents for linear cocycles. They proved the continuity of all Lyapunov exponents with respect to the fiber using exponential local uniform LDT estimates for linear cocycles in \cite{duarte2016lyapunov}. This is different from the one in our current paper which is the measure version of ACT. In \cite{duarte2020large}, they proved that the extremal Lyapunov exponents vary continuously with random matrix products for $\mathrm{GL}(2,\mathbb{R})$-finite valued cocycles. Moreover, if the simplicity of the Lyapunov exponents is assumed, their local H\"older continuity can be obtained. Lately, Duarte and Freijo \cite{duartefreijo2024continuity} established exponential uniform LDT estimates and H\"older continuity of the Lyapunov exponents with respect to fiber for random non-invertible cocycles with constant rank. Under the assumptions of quasi-irreducibility and the gap between the first two Lyapunov exponents, Baraviera and Duarte \cite{baraviera2019approximating} demonstrated that for random $\mathrm{GL}(d,\mathbb{R})$-valued cocycles, the first Lyapunov exponent is Lipschitz continuous on the fiber $A$ and the measure $\mu$ in the sense of the total variation norm.
\par In fact, the results outlined before were all established with the premise of compact support. In the case of possibly non-compact, Sánchez and Viana \cite{sanchez2018lyapunov} demonstrated that the first Lyapunov exponent of $\mathrm{SL}(2,\mathbb{R})$-valued cocycles over Bernoulli shifts is upper semi-continuous depending on measures in the Wasserstein distance, but not continuous for Wasserstein topology.
\par Extending a bit the scope of models, mixed random-quasiperiodic cocycles are natural generalization of product of random matrices. Set the model of mixed random-quasiperiodic $\mathrm{SL}(d,\mathbb{R})$-cocycles as a foundation, Cai, Duarte, and Klein \cite{cai2022mixed} established uniform base-LDT estimates in measure, and proved the uniform upper semicontinuity of the first Lyapunov exponent. Subsequently, for $\mathrm{SL}(d,\mathbb{R})$-cocycles over Bernoulli shifts in \cite{cai2023furstenberg}, they provided a Furstenberg-type formula to represent the first Lyapunov exponent and pointed out that it is continuous as the measure varies under the assumption of measure irreducibility, which can be deduced essentially from Kifer \cite{kifer2012ergodic}. Related to this, it was proved by Bezerra and Poletti \cite{bezerra2019random} that for finitely supported random quasi-periodic $\mathrm{SL}(2,\mathbb{R})$-valued cocycles, there exists an open and dense set of cocycles that are points of continuity for the Lyapunov exponents in the $C^0$ topology. Similar results are obtained in higher dimensions too. 
\par Compared with the results introduced above, our work transfers the difficulty of proving continuity of Lyapunov exponents to proving uniform LDT estimates. For example, Theorem \ref{thm1} can be directly applied to the continuity problem of Lyapunov exponents for product of random matrices.\par 
To formalize this connection, we first recall the notion of quasi-irreducibility: a cocycle $(\mu,A)$ is called \textit{quasi-irreducible} if either there is no proper subspace $V \subset \mathbb{R}^d$ which is invariant under all matrices of the cocycle generated by $\mathrm{supp}(\mu)$, or the only invariant proper subspace realizes the first Lyapunov exponent.
\begin{theorem}\label{product}
	Assume that the measure $\mu$ is quasi-irreducible and $L_1(\mu)>L_2(\mu)$, then the first Lyapunov exponent of random matrix products is H\"older continuous with respect to measure in the Wasserstein distance. 
\end{theorem}
\begin{remark}
	Actually, Dur\~aes \cite{pinto2021holder} has successfully proved the same result for $\mathrm{SL}(2,\mathbb{R})$-valued cocycles as Theorem \ref{product} in his master thesis in a completely different way. By employing the crucial assumption of measure quasi-irreducibility and the positivity condition of the Lyapunov exponents, he skillfully utilized the key property of the uniqueness of the stationary measure and combined it with the Furstenberg-Ledrappier formula, ultimately concluding that the first Lyapunov exponent is H\"older continuous in measure. It is worth emphasizing that our methodological framework diverges fundamentally from Dur\~aes's approach. In our current paper, we derive Theorem \ref{product} in higher dimensions as a direct corollary of our main Theorem \ref{thm1} when applied to random matrix products. In addition, our Theorem \ref{thm1.2} allows us to prove joint H\"older continuity with extra effort in obtaining uniform LDT estimates (in measure and fiber) for product of random matrices. However, we leave the proof to the readers for the sake of simplicity.
\end{remark}
\begin{remark}
	Our Theorem \ref{thm1} \& \ref{thm1.2} are applicable to cocycles that essentially generated by shifts, including mixed random-quasiperiodic cocycles mentioned above. Moreover, it is promisingly available for every deterministic cocycle as each deterministic one can be trivially rewritten as cocycles generated by shift with Dirac measure. However, the challenge now is converted into proving uniform base and fiber LDT estimates (in measure and fiber).
\end{remark}
	An important application of the result above arises in the study of random Schr\"odinger operators, where we can characterize the regularity of Lyapunov exponents as follows:
	\begin{corollary}\label{cor}
		Consider the random Schr\"odinger cocycle driven by the measure $\mu_E:=\int_{\mathbb{R}}\delta_{S_E(\omega)}d\rho(\omega)$ where $\rho\in \mathrm{Prob}_c(\mathbb{R})$. If $\rho$ is not a single Dirac, then the first Lyapunov exponent is H\"older continuous with respect to measure in Wasserstein distance.  
	\end{corollary}
	\begin{remark}
		When the energy parameter $E$ is slightly perturbed, the associated measure $\mu_E$ varies continuously in Wasserstein distance. This corollary naturally leads to the H\"older continuity of the Lyapunov exponent as a function of the energy. This is equivalent to the Hölder continuity of the Integrated Density of States (IDS), as established for example by Le Page \cite{le1989regularite} for one-parameter families of random matrices. And the quantitative bounds on the moduli of the continuity of the IDS were proved by Hislop-Marx \cite{hislop2020dependence} and Shamis \cite{shamis2021continuity}. In general, for most dynamical systems where a fiber-dependent measure can be constructed, our Theorem \ref{thm1} implies the fiber version of  ACT \cite{duarte2016lyapunov}.
	\end{remark}
\par The organization of this paper is as follows. In Section \ref{sec:pre} we explain some notations, introduce the definition and properties of the Wasserstein distance, and also state the general form of the Avalanche Principle. In Section \ref{sec:step}, we apply the Avalanche Principle to iterate the system. In each iteration, we utilize the Urysohn type lemma and the properties of the Wasserstein distance to provide reasonable estimates of the measure of the deviation set when measures vary on the base. Thereafter, the proof of our Theorem \ref{thm1} \& \ref{thm1.2} will be finalized in Section \ref{sec:continuity}. In Section \ref{apply} we directly apply our main theorem to random matrix products under necessary assumptions, obtaining the H\"older continuity of the first Lyapunov exponent.
\section{Preliminaries}\label{sec:pre}
\subsection{Lyapunov exponents and Wasserstein distance}
\par A cocycle $(\mu,A)$ is said to be integrable if $$\int_X\log^+\|A(x)\|d\mu(x)<+\infty.$$ Since $A$ is bounded, we have $\log^+\|A\|\in L^1$, hence by Kingman's subadditive ergodic theorem \cite{arnold1995random}, the first (or maximal) Lyapunov exponent of the cocycle exists, and it is denoted by $L_1(\mu,A)$, thus
\begin{align*}
    L_1(\mu,A)&=\lim_{n\to +\infty}\frac{1}{n}\log\|A^{(n)}(x)\|\quad \text{for $\mu$-a.e. $x\in X$}\\
    &=\lim_{n\to+\infty}\int_X\frac{1}{n}\log\|A^{(n)}(x)\|d\mu(x).
\end{align*}
Moreover, we call $L_1^{(n)}(\mu,A):=\displaystyle\int_X\frac{1}{n}\log\|A^{(n)}(x)\|d\mu(x)$ finite scale Lyapunov exponents. 
\par Additionally, using exterior powers of the cocycle $(\mu,A)$, the remaining Lyapunov exponents $L_j(\mu,A)$ can be expressed in terms of similar limits. For $1\leq k\leq j\leq d$ and $\mu$-a.e. $x\in X$,
\begin{align}\label{k1}
    \sum_{k=1}^j L_i(\mu,A)=\Lambda_j(\mu,A)=\lim_{n\to +\infty}\frac{1}{n}\log\|\wedge_jA^{(n)}(x)\|=\lim_{n\to +\infty}\frac{1}{n}\sum_{k=1}^j\log s_k(A^{(n)}(x)),
\end{align}
where $\wedge_jg$ denotes the $j$th exterior power of $g\in \mathrm{Mat}(d,\mathbb{R})$ and $\Lambda_0(g)=0$. Moreover, for $1\leq k\leq d$,
\begin{align}
	&s_1(\wedge_kg)=s_1(g)s_2(g)\cdots s_{k-1}(g)s_k(g),\label{k2}\\ 
	&s_2(\wedge_kg)=s_1(g)s_2(g)\cdots s_{k-1}(g)s_{k+1}(g). \label{k3}
\end{align}
\par To study the continuity of the Lyapunov exponents, we perturb the measure  within as small a neighborhood of $\mu$ as possible. Let $(X,\mu)$ and $(Y,\nu)$ be two probability spaces. A coupling of $\mu$ and $\nu$ is a measure $\pi$ on $X\times Y$ such that $\pi$ projects to $\mu$ and $\nu$ on the first and second coordinate, respectively. We endow the space $\mathscr{C}$ with the Wasserstein distance: for two probability measures $\mu$ and $\nu$ on $X$,
$$W_1(\mu,\nu)=\inf_{\pi\in\Pi(\mu,\nu)}\int_Xd(x,y)d\pi(x,y),$$
where $\Pi(\mu,\nu)$ is the set of all couplings of $\mu$ and $\nu$. 
\par In particular, Theorem 5.10(i) and Particular Case 5.4 in Villani \cite{villani2009optimal} together lead
to the useful duality formula for the Kantorovich–Rubinstein
distance: for any $\mu,\nu\in P(X),\xi\in \mathcal{C}^\alpha(X)$,
\begin{equation}\label{w1}
    W_1(\mu,\nu)=\sup_{\|\xi\|_{Lip}\leq 1}\left\{\int_X\xi d\mu-\int_X\xi d\nu\right\}.
\end{equation}
The distance $W_1(\mu,\nu)$ has Lipschitz continuity proved in the Theorem 7.29 of Villani\cite{villani2009optimal}:
\begin{equation}\label{w2}
	\left|\int_X\xi d\mu-\int_X\xi d\nu\right|<\|\xi\|_{Lip} W_1(\mu,\nu).
\end{equation}
\par It is well known that this distance metrizes the weak topology on $P(X)$. For further relevant details, we refer the readers to Chapter 6 of Villani  \cite{villani2009optimal}. From now on, we assume by default that the measure $\nu\in P(X)$ mentioned is an ergodic probability measure.
\par Below we state the Urysohn type lemma that will be used constantly in each iteration step to ensure the measure of the exceptional set does not vary too much when the base measure varies in $W_1$ topology.
\begin{lemma}\label{ury1}\cite{cai2022mixed}
	Let $X$ be a metric space and let $\mu$ be a Borel probability measure in $X$. Given a closed set $L\subset M$ and $\varepsilon>0$ there are an open set $D\supset L$ such that $\mu(D)<\mu(L)+\varepsilon$ and a Lipschitz continuous function $g:M\to[0,1]$ such that $\mathds{1}_L\leq g\leq \mathds{1}_D$.
\end{lemma}

\subsection{The Avalanche Principle}
\par For the purpose of proving the continuity of Lyapunov exponents for cocycles, we need to perform inductive processing across system iterations. The Avalanche Principle (AP), which correlates the norm growth of matrix products with the product of norms of individual matrices, provides the technical support for this goal.
\par The AP was initially introduced as a technical tool by Goldstein and Schlag in their study of matrix properties on $\mathrm{SL}(2, \mathbb{C})$ in \cite{goldstein2001holder}. After that, Schlag \cite{schlag2013regularity} generalized the AP to invertible matrices in $\mathrm{GL}(d, \mathbb{C})$. Then Duarte and Klein \cite{duarte2014continuity} extended the conclusions to arbitrary dimensions, providing a more precise theoretical statement for higher dimensions. Moreover, they offer a more general version in \cite{duarte2016lyapunov}, which is applicable to $\mathrm{Mat}(d,\mathbb{R})$ including non-invertible matrices.
\par We briefly review the content of the Avalanche Principle statement.
\par The gap ratio of the first two singular values of a matrix $g\in \mathrm{Mat}(d,\mathbb{R})$ is defined to be $$gr(g):=\frac{s_1(g)}{s_2(g)}\geq 1.$$
Moreover, $s_1(g)\cdot\cdots\cdot s_n(g)=\|\wedge_ng\|$.
\par Given a chain of linear mappings $\{g_j:V_j\to V_{j+1}\}_{0\leq j\leq n-1}$ where $V_i(0\leq i\leq n)$ represents Euclidean spaces of the same dimension $d$. We denote the composition of the first $i$ maps by $g^{(i)} := g_{i-1}\cdots g_1 g_0$. 
Throughout this paper, $a\lesssim b$ will stand for $a\leq Cb$ for some absolute constant $C$, and $a_n\asymp b_n$ means $0<\lim\limits_{n\to\infty}a_n/b_n<\infty$.
\begin{proposition}\label{ap}\cite{duarte2016lyapunov}
    There exists $c>0$ such that given $0<\varepsilon<1,0<\varkappa\leq c\varepsilon^2$ and $g_0,g_1,\cdots,g_{n-1}\in \mathrm{Mat}(d,\mathbb{R})$, if
    \begin{align}
        & \text{(gaps)}\quad gr(g_i)>\frac{1}{\varkappa} \qquad\qquad\text{for all}\quad 0\leq i\leq n-1,\label{ap1}\\
        &\text{(angles)}\ \frac{\|g_ig_{i-1}\|}{\|g_i\|\|g_{i-1}\|}>\varepsilon\qquad\text{for all}\quad 0\leq i\leq n-1.\label{ap2}
    \end{align}
    then
    \begin{align}
       \left|\log\|g^{(n)}\|+\sum_{i=1}^{n-2}\log\|g_i\|-\sum_{i=1}^{n-1}\log\|g_ig_{i-1}\|\right|\lesssim n\frac{\varkappa}{\varepsilon^2}.\label{ap4}
    \end{align}
\end{proposition}

\section{The Inductive Step and Other Technicalities}\label{sec:step}
The continuity of Lyapunov exponents cannot be easily proven directly. In this section, we introduce some preparatory work and important theoretical techniques for proving our main theorem. We will use uniform LDT estimates to prove the uniform continuity of the finite scale Lyapunov exponents. It ensures that the conditions of the AP are satisfied and allows us to establish inductive lemmas. These steps will enable the argument for the continuity theorem.
\subsection{Upper Semicontinuity of the First Lyapunov Exponent}
With the assumptions in this paper, it is evident that the first Lyapunov exponent is upper semi-continuous. Due to the necessity for proving the main theorem, we will present a stronger conclusion: the upper semicontinuity of the first Lyapunov exponent holds uniformly with respect to measure for a large enough set of phases, provided that the cocycle satisfies the uniform base-LDT estimates.
\begin{proposition}[nearly uniform upper semicontinuity]\label{prop1}
    Let $(\mu,A)\in\mathscr{C}$ be a continuous cocycle such that every observable $\xi\in \mathcal{C}^\alpha(X)(0<\alpha\leq 1)$ satisfies a uniform base-LDT estimate.\\
    (i) Assume that $L_1(\mu)>-\infty$.\\For every $\varepsilon>0$, there are $\delta=\delta(\mu,\varepsilon)>0,m=m(\mu,\varepsilon)\in\mathbb{N}$ and $c=c(\mu,\varepsilon)>0$, such that if $\nu\in P(X)$ with $W_1(\mu,\nu)<\delta$, and if $n\geq m$, then $$\frac{1}{n}\log\|A^{(n)}(x)\|\leq L_1(\mu)+\varepsilon$$ holds for all $x$ outside of a set of $\nu$-measure $<r(n)=e^{-cn}$.\\
    (ii) Assume that $L_1(\mu)=-\infty$.\\For every $t>0$, there are $\delta=\delta(\mu,t)>0,m=m(\mu,t)\in\mathbb{N}$ and $c=c(\mu,t)>0$ such that if $\nu\in P(X)$ with $W_1(\mu,\nu)<\delta$, and if $n\geq m$, then $$\frac{1}{n}\log\|A^{(n)}(x)\|\leq -t$$ holds for all $x$ outside of a set of $\nu$-measure $<r(n)$.
\end{proposition}
\begin{proof}
    Consider $C$ as the bound of $A$, i.e. $\|A(x)\|<C$. It deserves to be stated that $C$ represents a sufficiently large, positive and finite constant, which depends only on $\mu$ and varies slightly in different estimates during the proof process.
    \par (i) Fix $\varepsilon>0$. It's obvious that the sequence of functions $\log\|A^{(n)}(x)\|$ is sub-additive. By Kingman's subadditive ergodic theorem, $$\lim\limits_{n\to\infty}\frac{1}{n}\log\|A^{(n)}(x)\|=L_1(\mu) 
    \text{\quad for }\mu \text{ a.e. }x\in X,$$ hence $\text{for }\mu \text{ a.e. }x\in X$ we can define
    \begin{equation}
        n(x):=\min\{n\geq 1:\frac{1}{n}\log\|A^{(n)}(x)\|<L_1(\mu)+\varepsilon\}.
    \end{equation}
    \par Given $q\in\mathbb{N}$, let$$\mathcal{U}_q:=\{x:n(x)\leq q\}=\bigcup_{n=1}^q\{x:\frac{1}{n}\log\|A^{(n)}(x)\|<L_1(\mu)+\varepsilon\}.$$
    Since $A$ is continuous, the set $\mathcal{U}_q$ is open. Moreover, $\mathcal{U}_q\subset \mathcal{U}_{q+1}$ and $\cup_q\mathcal{U}_q$ has full $\mu$-measure. Therefore, there is $N=N(\mu,\varepsilon)$ such that $\mu(\mathcal{U}_N^\complement)<\varepsilon$.
    \par We fix this integer $N$ for the rest of the proof. If $x\in\mathcal{U}_N$ then $1\leq n(x)\leq N$ and 
    \begin{equation}
        \log\|A^{(n(x))}(x)\|\leq n(x)L_1(\mu)+n(x)\varepsilon.
    \end{equation}
    \par Let $m=m(\varepsilon,\mu):=\frac{CN}{\varepsilon}$.
    \par Fix $x\in X$ and define inductively for all $k\geq 1$ the sequence of phases $x_k=x_k(x)\in X$ and the sequence of integers $n_k=n_k(x)\in\mathbb{N}$ as follows:
    \begin{equation*}
    \begin{aligned}
       & x_1=x\qquad\qquad\qquad\quad n_1=
        \begin{cases}
         n(x_1)& \text{if $ x_1 \in \mathcal{U}_N $ } \\
         1& \text{if $ x_1 \notin \mathcal{U}_N $ }
        \end{cases}\\
        &x_2=T^{n_1}x_1\qquad\qquad\quad n_2=
        \begin{cases}
         n(x_2)& \text{if $ x_2 \in \mathcal{U}_N $ } \\
         1& \text{if $ x_2 \notin \mathcal{U}_N $ }
        \end{cases}\\
        &\cdots\\
       & x_{k+1}=T^{n_k}x_k\qquad\qquad n_{k+1}=
        \begin{cases}
         n(x_{k+1})& \text{if $ x_{k+1} \in \mathcal{U}_N $ } \\
         1& \text{if $ x_{k+1} \notin \mathcal{U}_N $ }.
        \end{cases}
    \end{aligned}
\end{equation*}
\par Note that for all $k\geq 1,x_{k+1}=T^{n_k+\cdots+n_1}x$ and $1\leq n_k\leq N$.
\par For any $n\geq m>N\geq n_1$, there is $p\geq 1$ such that $$n_1+\cdots+n_p\leq n< n_1+\cdots+n_p+n_{p+1},$$ so $n-n_{p+1}<n_1+\cdots+n_p\leq n$ and $n=n_1+\cdots+n_p+l$, where $0\leq l<n_{p+1}\leq N$.
\par Let $a_n(x):=\log\|A^{(n)}(x)\|$. Then clearly we have $a_n(x)\leq nC$,  and for all $n,l\geq 1$ and for every $x\in X$, $a_n(x)$ is a sub-additive process:$$a_{n+l}(x)\leq a_n(x)+a_l(T^nx).$$ This fact can be obtained from the Birkhoff ergodic theorem \cite{viana2016foundations} and compatibility of matrix norms.
\par By sub-additivity, we have
\begin{equation}
    \log\|A^{(n)}(x)\|=a_n(x)=a_{n_1+\cdots+n_p+l}(x)
    \leq \sum\limits_{k=1}^pa_{n_k}(x)+a_l(x_{p+1}).
\end{equation}
Obviously,
\begin{equation}
    a_l(x_{p+1})\leq lC<NC.
    \label{a1}
\end{equation}
\par For the earlier part, for every $1\leq k\leq p$, either $x_k\in\mathcal{U}_N$, so $n_k=n(x_k)$ and$$a_{n_k}(x_k)=\log\|A^{(n(x_k))}(x_k)\|\leq n(x_k)L_1(\mu)+n(x_k)\varepsilon,$$ or $x_k\notin\mathcal{U}_N$, so $n_k=1$ and $a_{n_k}(x_k)=\log\|A(x_k)\|\leq C.$
\par Therefore, for all $1\leq k\leq p$, 
\begin{align}
     a_{n_k}(x_k)&=a_{n_k}(x_k)\mathds{1}_{\mathcal{U}_N}(x_k)+a_{n_k}(x_k)\mathds{1}_{\mathcal{U}_N^\complement}(x_k)\nonumber\\
        &\leq (n(x_k)L_1(\mu)+n(x_k)\varepsilon)\mathds{1}_{\mathcal{U}_N}(x_k)+C\cdot\mathds{1}_{\mathcal{U}_N^\complement}(x_k)\nonumber\\
        &=(L_1(\mu)+\varepsilon)n_k+(C-L_1(\mu)-\varepsilon)\mathds{1}_{\mathcal{U}_N^\complement}(x_k)\nonumber\\
        &< (L_1(\mu)+\varepsilon)n_k+2C\cdot\mathds{1}_{\mathcal{U}_N^\complement}(x_k).
        \label{a2}
\end{align}
\par Then add up \eqref{a1} and \eqref{a2} for all $1\leq k\leq p$, we have:
\begin{equation*}
    \begin{aligned}
        \log\|A^{(n)}(x)\|&\leq (L_1(\mu)+\varepsilon)(n_1\cdots+n_p)+2C\cdot\sum_{k=1}^p\mathds{1}_{\mathcal{U}_N^\complement}(x_k)+CN\\
        &\leq n(L_1(\mu)+\varepsilon)+2C\cdot\sum_{j=0}^{n-1}\mathds{1}_{\mathcal{U}_N^\complement}(T^jx)+CN.
    \end{aligned}
\end{equation*}
Hence we can obtain the following: for all $x\in X$ and for all $n\geq m$, 
\begin{equation}
    \frac{1}{n}\log\|A^{(n)}(x)\|\leq L_1(\mu)+2\varepsilon+2C\frac{1}{n}\sum_{j=0}^{n-1}\mathds{1}_{\mathcal{U}_N^\complement}(T^jx).
\end{equation}
\par We use the Urysohn type lemma \ref{ury1}: for the closed set $\mathcal{U}_N^\complement$, there are an open set $D\supset\mathcal{U}_N^\complement$ such that $\mu(D)<\mu(\mathcal{U}_N^\complement)+\varepsilon$ and a Lipschitz continuous function $\xi:X\to [0,1]$ with $\|\xi\|_{Lip}=\frac{1}{\varepsilon}$ such that $\mathds{1}_{\mathcal{U}_N^\complement}\leq \xi\leq \mathds{1}_D$.
\par Applying the property of $W_1(\mu,\nu)$ and the uniform base-LDT to $\xi$, let $c\leq c'$, there are $\delta\leq\bar{\delta}, m\geq\bar{m}(\mu,\varepsilon)$ such that for all $n\geq m$, choosing $0<\delta\leq\varepsilon^2$, we get
\begin{equation*}
    \begin{aligned}
        \frac{1}{n}\sum_{j=0}^{n-1}\mathds{1}_{\mathcal{U}_N^\complement}(T^jx)&\leq  \frac{1}{n}\sum_{j=0}^{n-1}\xi(T^jx)\\&\leq \int_X\xi d\nu(x)+\varepsilon< \int_X\xi d\mu(x)+\frac{W_1(\mu,\nu)}{\varepsilon}+\varepsilon\\
        &<\int_X\mathds{1}_Dd\mu(x)+\frac{\delta}{\varepsilon}+\varepsilon=\mu(D)+\frac{\delta}{\varepsilon}+\varepsilon\\&<\mu(\mathcal{U}_N^\complement)+\varepsilon+\frac{\delta}{\varepsilon}+\varepsilon<4\varepsilon.
    \end{aligned}
\end{equation*}
holds for all $x$ outside of a set of $\nu$-measure $<r(n)$.
\par This ends the proof in the case $L_1(\mu)>-\infty$.
\par (ii) The case $L_1(\mu)=-\infty$. Consider $t>C+1$ be large enough. Applying again Kingman's subadditive ergodic theorem,  we define similarly as before:
\begin{align}
   & n(x):=\min\{n\geq 1:\frac{1}{n}\log\|A^{(n)}(x)\|<-2t\},\text{\quad for }\mu \text{ a.e. }x\in X\\
   & \mathcal{U}_q:=\{x:n(x)\leq q\}=\bigcup_{n=1}^q\{x:\frac{1}{n}\log\|A^{(n)}(x)\|<-2t\}.\nonumber
\end{align}
    Fix $N=N(\mu,\varepsilon)\in\mathbb{N}$. Then $\mu(\mathcal{U}_N^\complement)<\varepsilon=\frac{1}{t}$. If $x\in \mathcal{U}_N$, then $1\leq n(x)\leq N$ and 
\begin{equation}
    \log\|A^{(n(x))}(x)\|<-2tn(x).
\end{equation}
\par The rest of the proof  is similar to the case $L_1(\mu)>-\infty$. By stopping time argument in the same manner, we have for all $n\geq m$,
$$\frac{1}{n}\log\|A^{(n)}(x)\|\leq -2t+C\frac{1}{n}\sum_{j=0}^{n-1}\mathds{1}_{\mathcal{U}_N^\complement}(T^jx)-\frac{CN}{n}.$$
Then  for all $x$ outside of a set of $\nu$-measure $<r(n)$, we have $$\frac{1}{n}\log\|A^{(n)}(x)\|\leq -t.$$
%\par The rest of the proof  is similar to the case $L_1(\mu)>-\infty$, and readers can skip it. However, for the sake of completeness, the following part is a brief explanation of the proof:
\end{proof}
\begin{remark}
    If $L_1(\mu)=-\infty$, it is clear that $L_1$ is also lower semi-continuous at $\mu$, and hence continuous. Since $L_1(\mu)\geq L_2(\mu)\geq\cdots\geq L_d(\mu)$, each Lyapunov exponent is continuous at $\mu$. Therefore, from now on, we assume $L_1(\mu)>-\infty$.
\end{remark}
\begin{remark}
	Since the measure on the base varies within a neighborhood of $\mu$ (w.r.t. the distance $W_1(\nu,\mu)<\delta)$, it is necessary to estimate the measure of the same set under distinct base measures by using Urysohn type lemma here and the continuity of the fiber action $A$ is essential to make the result applicable. This is different from the situation in the articles of Duarte-Klein \cite{duarte2016lyapunov} where the base is fixed.
\end{remark}

The following lemma mainly applies Proposition \ref{prop1} and provides a lower bound on the gap between the first two singular values of the iterations of a cocycle. It ensures the gap condition of the AP satisfied and plays an important role in inductive iteration process.
\par We mention that if $(\mu,A)\in\mathscr{C}$ is such that $L_1(\mu)>L_2(\mu)\geq -\infty$, then we denote the gap between the first two Lyapunov exponents with $\kappa(\mu)$, i.e. $\kappa(\mu):=L_1(\mu)-L_2(\mu)>0$ when $L_2(\mu)>-\infty$, while if $L_2(\mu)=-\infty,\kappa(\mu)$ is a fixed, finite and large enough constant. 
\begin{lemma}\label{lem1}
    Let $(\mu,A)\in\mathscr{C}$ be a cocycle for which $L_1(\mu)>L_2(\mu)$ and let $\varepsilon>0$. There are $\delta_0=\delta_0(\mu,\varepsilon)>0,m_0=m_0(\mu,\varepsilon)\in\mathbb{N}$ and $c_0=c_0(\mu,\varepsilon)>0$ such that for all $\nu\in P(X)$ with $W_1(\mu,\nu)<\delta_0$ and for all $n\geq m_0$, if 
    \begin{equation}\label{b1}
        |L_1^{(n)}(\nu)-L_1^{(n)}(\mu)|<\theta,
    \end{equation}
   where $0<\theta<\varepsilon$, then for all phases $x$ outside a set of $\nu$-measure $<r_0(n)=e^{-c_0n}$ we have
    \begin{equation}\label{b2}
        \frac{1}{n}\log gr(A^{(n)}(x))>\kappa(\mu)-2\theta-3\varepsilon.
    \end{equation}
    Moreover,
    \begin{equation}\label{b3}
        L_1^{(n)}(\nu)-L_2^{(n)}(\nu)>(\kappa(\mu)-2\theta-3\varepsilon)(1-r_0(n)).
    \end{equation}
\end{lemma}
\begin{proof}
    Fix $\varepsilon>0$. Since $L_1(\mu)>L_2(\mu)$, the cocycle $(\mu,A)$ satisfies a uniform fiber-LDT estimate. Let $r_0(n)\geq r(n)+r''(n)$, and choose $\delta_0\leq\delta,m_0\geq m $ to ensure that both the uniform fiber-LDT estimate and Proposition \ref{prop1} hold.
    \par Note that $L_1(\wedge_2A)=L_1(\mu)+L_2(\mu)$, hence $L_1(\wedge_2A)>-\infty$ iff $L_2(\mu)>-\infty$. Using the Proposition \ref{prop1} on $\wedge_2A$ is indeed an effective approach.
    \par For any matrix $g\in \mathrm{Mat}(d,\mathbb{R})$, we have
    $$gr(g)=\frac{s_1(g)}{s_2(g)}=\frac{\|g\|^2}{\|\wedge_2g\|}\in (1,\infty].$$
    Then substitute matrix $A^{(n)}(x)$ into the above formula and we get
    \begin{equation}\label{b4}
         \frac{1}{n}\log gr(A^{(n)}(x))= 2\frac{1}{n}\log \|A^{(n)}(x)\|-\frac{1}{n}\log \|\wedge_2A^{(n)}(x)\|.
    \end{equation}
    \par Based on the uniform fiber-LDT estimate, the Kingman's subadditive ergodic theorem and assumption \eqref{b1}, we can derive that
    \begin{equation}\label{b5}
        \frac{1}{n}\log \|A^{(n)}(x)\|>L_1^{(n)}(\nu)-\varepsilon>L_1^{(n)}(\mu)-\theta-\varepsilon\geq L_1(\mu)-\theta-\varepsilon
    \end{equation}
    hold for all $x$ outside a set $E_1$ of $\nu$-measure $<r''(n)$.
    \par Concerning the last term of \eqref{b4}, applying Proposition \ref{prop1} to $\wedge_2A$, we will obtain an upper bound on $\frac{1}{n}\log\|\wedge_2A^{(n)}(x)\|$.
    \par If $L_2(\mu)>-\infty$, so $L_1(\wedge_2A)>-\infty$, then for all $x$ outside a set $E_2$ of $\nu$-measure $<r(n)$, we have 
    \begin{equation}\label{b6}
        \frac{1}{n}\log\|\wedge_2A^{(n)}(x)\|<L_1(\mu)+L_2(\mu)+\varepsilon,
    \end{equation}
    hence, combine \eqref{b4}, \eqref{b5} and \eqref{b6},
    $$\frac{1}{n}\log gr(A^{(n)}(x))>\kappa(\mu)-2\theta-3\varepsilon$$
    holds for $x$ outside a set $E_1\cup E_2$ of $\nu$-measure $<r_0(n)$.
    \par Therefore, considering $E=E_1\cup E_2$, and integrating in $x$ w.r.t. $\nu$, we obtain
    \begin{equation*}
    \begin{aligned}
        L_1^{(n)}(\nu)- L_2^{(n)}(\nu)&=\int_X\frac{1}{n}\log gr(A^{(n)}(x))d\nu(x)\\&=\int_{E^\complement}\frac{1}{n}\log gr(A^{(n)}(x))d\nu(x)+\int_{E}\frac{1}{n}\log gr(A^{(n)}(x))d\nu(x)\\&>(\kappa(\mu)-2\theta-3\varepsilon)(1-r_0(n)),
    \end{aligned}
    \end{equation*}
    since $gr(A^{(n)}(x))\in (1,+\infty]$, we deduce clearly that $\displaystyle\int_{E}\frac{1}{n}\log gr(A^{(n)}(x))d\nu(x)>0$, then \eqref{b3} is proven.
    \par If $L_2(\mu)=-\infty$, so $L_1(\wedge_2A)=-\infty$, let $t=t(\mu):=-2L_1(\mu)+\kappa(\mu)$. Then for all $x$ outside a set $E_3$ of $\nu$-measure $<r(n)$, we have 
    \begin{equation}\label{b7}
        \frac{1}{n}\log\|\wedge_2A^{(n)}(x)\|<L_1(\mu)+L_2(\mu)+\varepsilon,
    \end{equation}
    hence, combine \eqref{b4}, \eqref{b5} and \eqref{b7},
    $$\frac{1}{n}\log gr(A^{(n)}(x))>\kappa(\mu)-2\theta-3\varepsilon$$
    holds for $x$ outside a set $E_1\cup E_3$ of $\nu$-measure  $<r_0(n)$ .
    \par Moreover, we can similarly conclude \eqref{b3}.
\end{proof}
\subsection{Finite Scale Continuity}
We prove that the finite scale Lyapunov exponent is uniformly continuous with respect to measure when the scale is fixed, but it is not the same as continuity of the true Lyapunov exponent. In the following inductive steps, it can be demonstrated through iteration that the continuous behavior holds as the scale tends towards infinity.
\begin{proposition}[finite scale uniform continuity]\label{prop2}
     Let $(\mu,A)\in\mathscr{C}$ be a cocycle for which $L_1(\mu)>L_2(\mu)$. There are $\delta_1=\delta_1(\mu)>0,m_1=m_1(\mu)\in\mathbb{N},C_1=C_1(\mu)>0$ and $c_1=c_1(\mu)>0$ such that for any two measures $\nu_1,\nu_2\in P(X)$ with $W_1(\mu,\nu_i)<\delta_1$, where $i=1,2$, if $n\geq m_1$ and $W_1(\nu_1,\nu_2)<e^{-C_1n}$, then
     \begin{equation}\label{f1}
         |L_1^{(n)}(\nu_1)-L_1^{(n)}(\nu_2)|<r_1(n)=e^{-c_1n}.
     \end{equation}
\end{proposition}
\begin{proof}
   Fix $\varepsilon:=\kappa(\mu)/10>0$.
    \par Since $L_1(\mu)>L_2(\mu)$, choose appropriate $\delta_1\leq\delta_0,m_1\geq m_0,r_1(n)\geq 2r''(n)$ such that we can apply the uniform fiber-LDT estimate to complete our proof. Pick $C_1$ such that if $m_1>\bar{m},\ e^{-C_1m_1}<\delta_1$ which ensures the conditions of Lemma \ref{lem1} satisfied.
    \par We assume that the bound of $A$ is $e^{C}$ where $C=C(\mu)>0$. Relating to $C_0$ norm, for all $n\geq 1$,
    $$|L_1^{(n)}(\nu)|\leq\left\|\frac{1}{n}\log\|A^{(n)}(x)\|\right\|_{C_0}\leq C.$$
    \par Fix arbitrary $n\geq m_1$ and $\nu_i$ with $W_1(\nu_i,\mu)\leq\delta_1$. Apply the uniform fiber-LDT estimate and conclude that for all $x$ outside a set $E_i\ (i=1,2)$ with $\nu_i(E_i)<r''(n)$, 
    \begin{equation}\label{f2}
        \frac{1}{n}\log\|A^{(n)}(x)\|<L_1^{(n)}(\nu_i)+\varepsilon\leq C+\varepsilon.
    \end{equation}
    \par Let $E=E_1\cup E_2$, so $\nu_i(E)<r_1(n)$ where $i=1,2$.
    \par In fact, for instance, $\nu_1(E)\leq\nu_1(E_2)+\nu_1(E_1)<\nu_1(E_2)+r''(n)$. Since $A$ is continuous, $E_2=\{x\in X:|\frac{1}{n}\log\|A^{(n)}(x)\|-L_1^{(n)}(\nu_2)|>\varepsilon\}$ is an open set clearly, then $E_2^\complement$ is a closed set. 
    \par Therefore we can use the Urysohn Type lemma \ref{ury1}, there are an open set $D\supset E_2^\complement$ such that $\nu_1(D)<\nu_1(E_2^\complement)+\tau$ where $\tau>0$ and a Lipschitz continuous function $\xi:X\to [0,1]$ with $\|\xi\|_{Lip}=\frac{1}{\tau}$ such that $\mathds{1}_{E_2^\complement}\leq \xi\leq \mathds{1}_D$. 
    \par Hence
    \begin{equation*}
    \begin{aligned}
        \nu_2(E_2^\complement)&\leq \int_X\xi d\nu_2(x)< \int_X\xi d\nu_1(x)+\frac{W_1(\nu_1,\nu_2)}{\tau}\\
        &<\nu_1(D)+\frac{e^{-C_1n}}{\tau}<\nu_1(E_2^\complement)+\tau+\frac{e^{-C_1n}}{\tau}.
    \end{aligned}
\end{equation*}
Let $\tau=e^{\frac{-C_1n}{2}}$ and choose $C_1>2c_1$ large enough so that $\frac{e^{-C_1n}}{\tau}+\tau<r_1(n)$, then
\begin{equation}\label{f3}
    \nu_2(E_2^\complement)-\nu_1(E_2^\complement)<r_1(n).
\end{equation}
Then $\nu_1(E)<2r_1(n)$, similarly, $\nu_2(E)<2r_1(n)$.
\par If $x\in E$, according to basic inequality, we have 
\begin{equation}\label{f4}
    \begin{aligned}
        &\left|\int_{E}\frac{1}{n}\log\|A^{(n)}(x)\|d\nu_1(x)-\int_{E}\frac{1}{n}\log\|A^{(n)}(x)\|d\nu_2(x)\right|\\&\leq \int_{E}\left|\frac{1}{n}\log\|A^{(n)}(x)\|\right|d\nu_1(x)+\int_{E}\left|\frac{1}{n}\log\|A^{(n)}(x)\|\right|d\nu_2(x)\\
        &\leq\left\|\frac{1}{n}\log\|A^{(n)}(x)\|\right\|_{C_0}\cdot\nu_1(E)+\left\|\frac{1}{n}\log\|A^{(n)}(x)\|\right\|_{C_0}\cdot\nu_2(E)\\
        &\leq 4C\cdot r_1(n).
    \end{aligned}
\end{equation}
\par If $x\in E^\complement$, combine \eqref{f2} - \eqref{f3}, we obtain
\begin{equation}\label{f5}
    \begin{aligned}
        &\left|\int_{E^\complement}\frac{1}{n}\log\|A^{(n)}(x)\|d\nu_1(x)-\int_{E^\complement}\frac{1}{n}\log\|A^{(n)}(x)\|d\nu_2(x)\right|\\&\leq (C+\varepsilon)\cdot\left|\int_{E^\complement}d\nu_1(x)-\int_{E^\complement}d\nu_2(x)\right|\\
        &\leq(C+\varepsilon)\cdot(\frac{e^{-C_1n}}{\tau}+\tau)
        \leq (C+\varepsilon)r_1(n).
    \end{aligned}
\end{equation}
\par Therefore, sum up \eqref{f4} and \eqref{f5}, 
$$|L_1^{(n)}(\nu_1)-L_1^{(n)}(\nu_2)|<r_1(n).$$
\end{proof}
\subsection{The Inductive Step Procedure}
We will establish the main tool for proving our ACT here, i.e. the inductive step procedure which relies on the Avalanche Principle \ref{ap} and uniform fiber-LDT estimates \ref{de4}. We first present two lemmas, which have been proven in Duarte-Klein \cite{duarte2016lyapunov}.
\begin{lemma}\label{lem2}\cite{duarte2016lyapunov}
    Let $(\nu,A)\in\mathscr{C}$ be a cocycle. If $n_0,n_1,n,s\in \mathbb{N}$ are such that $n_1=n\cdot n_0+s$ and $0\leq s\leq n_0$, then 
    \begin{equation}
        -2C\frac{n_0}{n_1}+L_1^{((n+1)n_0)}(\nu)\leq L_1^{(n_1)}(\nu)\leq L_1^{(nn_0)}(\nu)+2C\frac{n_0}{n_1}.
    \end{equation}
\end{lemma}
\begin{lemma}\label{lem3}\cite{duarte2016lyapunov}
    Let $(\nu,A)\in\mathscr{C}$ for all $\nu\in P(X)$ satisfies a fiber-LDT estimate. Let $q_1,q_2,n\in\mathbb{N},c=c(\nu)>0$ and $\eta>0$ be such that $q_i\geq n\geq \bar{n}$ where $i=1,2$, and$$|L_1^{(q_1+q_2)}(\nu)-L_1^{(q_i)}(\nu)|<\eta.$$Then
    \begin{equation}
        \frac{\|A^{(q_1+q_2)}(x)\|}{\|A^{(q_2)}(T^{q_1}x)\|\|A^{(q_1)}(x)\|}>e^{-(q_1+q_2)(\eta+2\varepsilon)}
    \end{equation}
    holds for all $x$ outside a set of $\nu$-measure $<3r(n)$.
\end{lemma}
\par According to Lemma \ref{lem2}, for two consecutive scales $n_0$ and $n_1$, if $n$ is a multiple of $n_0$, these two finite scale Lyapunov exponents differ by a small value, approximately of order $\frac{n_0}{n_1}$. Moreover, the Lemma \ref{lem3} is used to ensure that the angle condition for the avalanche principle is satisfied.
\begin{proposition}[inductive step procedure]\label{prop3}
    Let $(\mu,A)\in\mathscr{C}$ such that $L_1(\mu)>L_2(\mu)$. Fix $\varepsilon>0$, there are $C_2=C_2(\mu)>0,\delta_2=\delta_2(\mu,\varepsilon),m_2=m_2(\mu,\varepsilon)\in\mathbb{N}\text{ and }c_2=c_2(\mu,\varepsilon)>0$ such that for any $n_0\geq m_2$, if the inequalities
    \begin{align}
         &(a)\quad L_1^{(n_0)}(\nu)-L_1^{(2n_0)}(\nu)<\eta_0\label{i1}\\
         &(b)\quad |L_1^{(n_0)}(\nu)-L_1^{(n_0)}(\mu)|<\theta_0\label{i2}
    \end{align}
    hold for a measure $\nu\in P(X)$ with $W_1(\mu,\nu)<\delta_2$, and if the positive numbers $\eta_0,\theta_0$ satisfy
    \begin{equation}\label{i3}
            2\theta_0+4\eta_0<\kappa(\mu)-12\varepsilon,
    \end{equation}
    then for any integer $n_1$ such that 
    \begin{equation}\label{i4}
            n_0^{1+}\leq n_1\leq\psi(n_0),
    \end{equation}
    where $\psi(n_0)=n_0\cdot [r_2(n_0)]^{-(\frac{1}{2}-)}$, we have
    \begin{equation}\label{i5}
           |L_1^{(n_1)}(\nu)+L_1^{(n_0)}(\nu)-2L_1^{(2n_0)}(\nu)|<C_2\frac{n_0}{n_1}.
    \end{equation}
    Furthermore,
     \begin{align}
         &(a)\quad L_1^{(n_1)}(\nu)-L_1^{(2n_1)}(\nu)<\eta_1\label{i6}\\
         &(b)\quad |L_1^{(n_1)}(\nu)-L_1^{(n_1)}(\mu)|<\theta_1,\label{i7}
     \end{align}
    where
     \begin{align}
       & \theta_1=\theta_0+4\eta_0+C_2\frac{n_0}{n_1},\label{i8}\\
        &\eta_1=C_2\frac{n_0}{n_1}.\label{i9}
    \end{align}
\end{proposition}
\begin{proof}
    Throughout this proof, $C_2$ denotes a positive, finite constant, sufficiently large and depending only on $\mu$, whose exact value may vary between different estimates. 
    \par Since $L_1(\mu)>L_2(\mu)$,  we can adjust the parameters $\delta_2\leq\delta_0,m_2\geq m_0,c_2<c_0$ to ensure that both the uniform fiber-LDT estimates and Lemma \ref{lem1} are satisfied. Let $r_2(n)=e^{-c_2n}\geq r_0(n)+3r''(n)$.
    \par From \eqref{i1}, applying Lemma \ref{lem3}, we have:
    \begin{equation}\label{i10}
        \frac{\|A^{(2n_0)}(x)\|}{\|A^{(n_0)}(T^{n_0}x)\|\|A^{(n_0)}(x)\|}>e^{-2n_0(\eta_0+2\varepsilon_{n_0})}\geq e^{-n_0(2\eta_0+4\varepsilon)}:=\varepsilon_{ap}
    \end{equation}
    holds for all $x$ outside a set $E_{n_0}^1$ of $\nu$-measure $<3r^{''}(n_0)$.
    \par Moreover, due to the assumption \eqref{i2}, applying Lemma \ref{lem1}: for all $x$ outside a set $E_{n_0}^2$ of $\nu$-measure $<r_0(n_0)$, we have
    $$\frac{1}{n}\log gr(A^{(n_0)}(x))>\kappa(\mu)-2\theta_0-3\varepsilon,$$
    then
    \begin{equation}\label{i11}
        gr(A^{(n_0)}(x))>e^{n_0(\kappa(\mu)-2\theta_0-3\varepsilon)}:=\frac{1}{\varkappa_{ap}}.
    \end{equation}
    These estimates \eqref{i10} and \eqref{i11} will ensure that the angle and gap conditions in the AP are satisfied.
    \par Let $E_{n_0}=E_{n_0}^1\cup E_{n_0}^2$,  so  $\nu(E_{n_0})<r_2(n_0)$.
    \par Moreover, let $n_1$ be an integer such that $n_0^{1+}\leq n_1\leq \psi(n_0)$, and choose $m_2$ large enough. We assume that $n_1=n\cdot n_0$ for some $n\in\mathbb{N}$.
    \par For every $0\leq i\leq n-1$, define
    $$g_i=g_i(x):=A^{(n_0)}(T^{in_0}x).$$
    Then we can get clearly $g^{(n)}=g_{n-1}\cdot\cdots\cdot g_0=A^{(n\cdot n_0)}(x)=A^{(n_1)}(x)$ and $g_ig_{i-1}=A^{(2n_0)}(T^{(i-1)n_0}x)$ for all $1\leq i\leq n-1$.
    \par Considering $\bar{E}_{n_0}:=\bigcup_{i=0}^{n-1}T^{-in_0}E_{n_0}$, so $\nu(\bar{E}_{n_0})<nr_2(n_0)$. If $x\notin\bar{E}_{n_0}$, then $T^{in_0}x\notin E_{n_0}$, so 
    \begin{align*}
        & gr(g_i)>\frac{1}{\varkappa_{ap}} \qquad\text{for all}\quad 0\leq i\leq n-1,\\
        &\frac{\|g_ig_{i-1}\|}{\|g_i\|\|g_{i-1}\|}>\varepsilon_{ap}\qquad\text{for all}\quad 0\leq i\leq n-1.
    \end{align*}
    Note that $\varkappa_{ap}\ll\varepsilon_{ap}^2$. In fact, we may obtain this inequality easily from the assumption \eqref{i3}. Therefore, applying the avalanche principle, we have
    $$\left|\log\|g^{(n)}\|+\sum_{i=1}^{n-2}\log\|g_i\|-\sum_{i=1}^{n-1}\log\|g_ig_{i-1}\|\right|\lesssim n_0\frac{\varkappa_{ap}}{\varepsilon_{ap}^2}.$$
    Note that $\frac{\varkappa_{ap}}{\varepsilon_{ap}^2}=e^{-n_0(\kappa(\mu)-4\eta_0-2\theta_0-11\varepsilon)}<e^{-\varepsilon n_0}$. Then $\frac{\varkappa_{ap}}{\varepsilon_{ap}^2}<r_2(n_0)$ since we choose $c_0>c_2>\varepsilon>0$.
    \par For all $x\notin\bar{E}_{n_0}$, where $\nu(\bar{E}_{n_0})<nr_2(n_0)$, using the AP, we get
    \begin{equation}\label{i12}
       \begin{split}
           \Biggl|\log\|A^{(n_1)}(x)\|&+\sum_{i=1}^{n-2}\log\|A^{(n_0)}(T^{in_0}x)\|\\&-\sum_{i=1}^{n-1}\log\|A^{(2n_0)}(T^{(i-1)n_0}x)\|\Biggr|\leq nr_2(n_0).
       \end{split}
    \end{equation}
    \par Divide both side of the above inequality by $n_1=n\cdot n_0$, for $x\notin\bar{E}_{n_0}$, we have
    \begin{equation}
       \begin{split}
           |f(x)|=&\Biggl|\frac{1}{n_1}\log\|A^{(n_1)}(x)\|+\frac{1}{n}\sum_{i=1}^{n-2}\frac{1}{n_0}\log\|A^{(n_0)}(T^{in_0}x)\|\\&-\frac{2}{n}\sum_{i=1}^{n-1}\frac{1}{2n_0}\log\|A^{(2n_0)}(T^{(i-1)n_0}x)\|\Biggr|\leq \frac{1}{n_0}r_2(n_0).
       \end{split}
    \end{equation}
   Then integrate in $x$ w.r.t. $\nu$,
   $$\int_Xf(x)d\nu(x)=L_1^{(n_1)}(\nu)+\frac{n-2}{n}L_1^{(n_0)}(\nu)-\frac{2(n-1)}{n}L_1^{(2n_0)}(\nu).$$
   \par Considering the $C_0$ norm and the assumption \eqref{i4}, we have
   \begin{equation*}
     \begin{aligned}
         \int_X|f(x)|d\nu(x)
&=\int_{\bar{E}_{n_0}^\complement}|f(x)|d\nu(x)+\int_{\bar{E}_{n_0}}|f(x)|d\nu(x)\\
        &\leq r_2(n_0)(1-nr_2(n_0))+\|f(x)\|_{C_0}\cdot\nu(\bar{E}_{n_0})\\
        &\leq r_2(n_0)+C_2\cdot\nu(\bar{E}_{n_0})
        =r_2(n_0)+C_2\cdot (nr_2(n_0))<C_2\frac{n_0}{n_1}.
       \end{aligned}
   \end{equation*}
   \par Therefore,
   \begin{equation*}
       \begin{aligned}
           \Biggl|\int_Xf(x)d\nu(x)\Biggr|&=|L_1^{(n_1)}(\nu)+\frac{n-2}{n}L_1^{(n_0)}(\nu)-\frac{2(n-1)}{n}L_1^{(2n_0)}(\nu)|\\&=|L_1^{(n_1)}(\nu)+L_1^{(n_0)}(\nu)-2L_1^{(2n_0)}(\nu)-\frac{2}{n}(L_1^{(n_0)}(\nu)-L_1^{(2n_0)}(\nu))|\\&\leq \int_X|f(x)|d\nu(x)<C_2\frac{n_0}{n_1},
       \end{aligned}
   \end{equation*}
   Further, the following conclusion can be reached:
   \begin{equation*}
       \begin{aligned}
           |L_1^{(n_1)}(\nu)+L_1^{(n_0)}(\nu)-2L_1^{(2n_0)}(\nu)|&<C_2\frac{n_0}{n_1}-\frac{2}{n}[L_1^{(n_0)}(\nu)-L_1^{(2n_0)}(\nu)]\\&<C_2\frac{n_0}{n_1}-\frac{2}{n}\cdot\eta_0\leq C_2\frac{n_0}{n_1}.
       \end{aligned}
   \end{equation*}
   \par Hence for the particular scales $n_1'=n\cdot n_0$ and $n_1''=(n+1)n_0$, we may obtain similar conclusions. Moreover, for any scale $n_1=n\cdot n_0+s$, where $0\leq s\leq n_0$, applying Lemma \ref{lem2}, we get
   \begin{equation}\label{i13}
       \begin{aligned}
          &|L_1^{(n_1)}(\nu)+L_1^{(n_0)}(\nu)-2L_1^{(2n_0)}(\nu)|\\&\leq |L_1^{(n_1)}(\nu)+L_1^{(n_0)}(\nu)-2L_1^{(2n_0)}(\nu)+2C_2\frac{n_0}{n_1}|\leq C_2\frac{n_0}{n_1}.
       \end{aligned}
   \end{equation}
   \par Considering $n_0^{1+}\leq n_1\leq\frac{1}{2}\psi(n_0)$, the above argument will hold for $2n_1$ instead $n_1$, then 
   \begin{equation*}
       \begin{aligned}
           L_1^{(n_1)}(\nu)-L_1^{(2n_1)}(\nu)&<L_1^{(n_1)}(\nu)+L_1^{(n_0)}(\nu)-2L_1^{(2n_0)}(\nu)+C\frac{n_0}{n_1}\\&< C_2\frac{n_0}{n_1}:= \eta_1,
       \end{aligned}
   \end{equation*}
   that is \eqref{i6} proven.
   \par We can rewrite \eqref{i13} as follows:
   $$|L_1^{(n_1)}(\nu)-L_1^{(n_0)}(\nu)+2[L_1^{(n_0)}(\nu)-L_1^{(2n_0)}(\nu)]|< C_2\frac{n_0}{n_1}.$$
   It also applies to $\mu$, then \eqref{i7} establishes,
   \begin{equation*}
       \begin{split}
           |L_1^{(n_1)}(\nu)-L_1^{(n_1)}(\mu)|=&\Biggl|L_1^{(n_1)}(\nu)-L_1^{(n_0)}(\nu)+2[L_1^{(n_0)}(\nu)-L_1^{(2n_0)}(\nu)]\\&-\big(L_1^{(n_1)}(\mu)-L_1^{(n_0)}(\mu)+2[L_1^{(n_0)}(\mu)-L_1^{(2n_0)}(\mu)]\big)\\&+L_1^{(n_0)}(\nu)-L_1^{(n_0)}(\mu)\\&+2[L_1^{(n_0)}(\nu)-L_1^{(2n_0)}(\nu)]-2[L_1^{(n_0)}(\mu)-L_1^{(2n_0)}(\mu)]\Biggr|\\&< C_2\frac{n_0}{n_1}+\theta_0+4\eta_0:=\theta_1.
       \end{split}
    \end{equation*}
\end{proof}
\section{Continuity Theorem}\label{sec:continuity}
Repeat the previous iterative process and we will obtain the continuity and modulus of continuity of Lyapunov exponents with respect to measures. In this section, we present the complete proof of our main theorem.
\subsection{General Continuity Theorem}\label{sec:gen}
By applying the finite scale uniform continuity Proposition \ref{prop2} and the inductive step Proposition \ref{prop3}, we will establish the continuity of the first Lyapunov exponent.

\begin{theorem}\label{thm2}
   Given an ergodic MPDS $(X,\mu,\mathcal{F},T)$, a space of continuous cocycle $\mathscr{C}$, we assume that:\\
   (1) every $\xi\in \mathcal{C}^\alpha(X)(0<\alpha\leq 1)$ satisfies a uniform base-LDT estimate.\\
   (2) every cocycle $(\mu,A)\in\mathscr{C}$ for which $L_1(\mu)>L_2(\mu)$ satisfies a uniform fiber-LDT estimate.
   \par Then there exist a neighborhood $\mathcal{V}_\mu$ of $\mu$ and $\delta>0$ such that for any $\nu\in\mathcal{V}_\mu$ with $W_1(\mu,\nu)<\delta$ the map $P(X)\ni\nu\mapsto L_1(\nu)$ is continuous.
\end{theorem} 
\begin{proof}
    Let $0<\varepsilon<\kappa(\mu)/50$ be arbitrary fixed.
    \par Since $L_1^{(n)}(\mu)\to L_1(\mu)$ as $n\to\infty$, there is $m_3=m_3(\mu,\varepsilon)\in\mathbb{N}$ such that for all $n\geq m_3$ we have 
    \begin{equation}\label{g1}
        L_1^{(n)}(\mu)-L_1^{(2n)}(\mu)<\varepsilon.
    \end{equation}
    \par Let $r_3(n)=e^{-c_3n}\geq e^{-c_1n}+e^{-c_2n}=r_1(n)+r_2(n)$ be the corresponding deviation measure function with $c_3<\min\{c_1,c_2\}$. Let $\delta\leq\min\{\delta_1,\delta_2\}$ and $C_1,C_2$ be the constants in Proposition \ref{prop2} and Proposition \ref{prop3} respectively.
    \par Let the scale $n_0\geq\max\{m_1,m_2,m_3\}$ . Moreover, assume $n_0$ to be large enough so that $e^{-C_12n_0}<\delta_0,r_3(n_0)<\varepsilon$ and $n^{0+}\ll [r_3(n)]^{-(\frac{1}{2}-)}$ for $n\geq n_0$ and $C_2\frac{n_0}{n_0^{1+}}<\delta$.
    \par Let $\delta:=e^{-C_12n_0}$ and let $\nu\in P(X)$ with $W_1(\mu,\nu)<\delta$.
    \par Since $\delta=e^{-C_12n_0}<e^{-C_1n_0}$, we can apply Proposition \ref{prop2} at scale $2n_0$ and $n_0$ and get
    \begin{align}
        &|L_1^{(n_0)}(\nu)-L_1^{(n_0)}(\mu)|<r_3(n_0):=\theta_0<\varepsilon \label{g2}\\
        &|L_1^{(2n_0)}(\nu)-L_1^{(2n_0)}(\mu)|<r_3(n_0)=\theta_0\label{g3}
    \end{align}
    \par Then
    \begin{equation}\label{g4}
        |L_1^{(n_0)}(\nu)-L_1^{(2n_0)}(\nu)|<2r_3(n_0)+\varepsilon:=\eta_0<3\varepsilon
    \end{equation}
    Moreover, $2\theta_0+4\eta_0<2\varepsilon+12\varepsilon=14\varepsilon<\kappa(\mu)-12\varepsilon$. Then for any integer $n_1$ satisfying $n_0^{1+}\leq n_1\leq \psi(n_0)$, we can apply the Proposition \ref{prop3} and conclude that 
     \begin{align}
         & L_1^{(n_1)}(\nu)-L_1^{(2n_1)}(\nu)<\eta_1\label{g6}\\
         &|L_1^{(n_1)}(\nu)-L_1^{(n_1)}(\mu)|<\theta_1,\label{g7}
     \end{align}
     where
     \begin{align}
       & \theta_1=\theta_0+4\eta_0+C_2\frac{n_0}{n_1},\\
        &\eta_1=C_2\frac{n_0}{n_1}.
    \end{align}
    \par Furthermore, $2\theta_1+4\eta_1=\theta_0+4\eta_0+6C_2\frac{n_0}{n_1}<26\varepsilon+\varepsilon<\kappa(\mu)-12\varepsilon$. Then for $n_1^{1+}\leq n_2\leq \psi(n_1)$, by Proposition \ref{prop3} again we get
    \begin{align*}
         & L_1^{(n_2)}(\nu)-L_1^{(2n_2)}(\nu)<\eta_2\\
         &|L_1^{(n_2)}(\nu)-L_1^{(n_2)}(\mu)|<\theta_2,
     \end{align*}
     where
     \begin{align*}
       & \theta_2=\theta_1+4\eta_1+C_2\frac{n_1}{n_2},\\
        &\eta_2=C_2\frac{n_1}{n_2}.
    \end{align*}
    Note also that $2\theta_2+4\eta_2<\kappa(\mu)-12\varepsilon$.
    \par Continue this procedure. We choose a scale $n_k^{1+}\leq n_{k+1}\leq \psi(n_k)$ and from Proposition \ref{prop3} we have:
    \begin{align}
         & L_1^{(n_{k+1})}(\nu)-L_1^{(2n_{k+1})}(\nu)<\eta_{k+1}\label{g9}\\
         &|L_1^{(n_{k+1})}(\nu)-L_1^{(n_{k+1})}(\mu)|<\theta_{k+1},\label{g10}
     \end{align}
     where
    \begin{equation}\label{g11}
        \eta_{k+1}=C_2\frac{n_k}{n_{k+1}},
    \end{equation}
    and
    \begin{align*}
       \theta_{k+1}&=\theta_k+4\eta_k+C_2\frac{n_k}{n_{k+1}}\\
       &=(\theta_0+4\eta_0)+5C_2\sum_{i=0}^{k-1}\frac{n_i}{n_{i+1}}+C_2\frac{n_k}{n_{k+1}}\\&<(\theta_0+4\eta_0)+5C_2\sum_{i=0}^{\infty}\frac{n_i}{n_{i+1}}\\&<(\theta_0+4\eta_0)+10C_2\frac{n_0}{n_0^{1+}}<(\varepsilon+12\varepsilon)+10\varepsilon=23\varepsilon.
    \end{align*}
    Hence,
    \begin{equation}\label{g12}
        \theta_{k+1}<23\varepsilon.
    \end{equation}
    Moreover, $2\theta_{k+1}+4\eta_{k+1}<\kappa(\mu)-12\varepsilon$, which guarantees that the inductive process can run indefinitely.
    \par Take the limit as $k\to\infty$ in \eqref{g10}, and we get
    $$|L_1(\nu)-L_1(\mu)|\leq 23\varepsilon,$$
    this implies the continuity of $L_1(\nu)$.
\end{proof}
\par As a matter of fact, if $(\mu,A)$ has a different gap pattern, for instance, $L_1(\mu)=L_2(\mu)>L_3(\mu)\geq\cdots L_d(\mu)$, consider instead the $\wedge_2A$, hence
\begin{align*}
	L_1(\wedge_2A)-L_2(\wedge_2A)&=L_1(\mu)+L_2(\mu)-(L_1(\mu)+L_3(\mu))\\
	&=L_2(\mu)-L_3(\mu)>0.
\end{align*} 
Applying Theorem \ref{thm2} to $\wedge_2A$, we conclude that $L_1(\nu)+L_2(\nu)$ is continuous. Either $L_1(\nu)>L_2(\nu)$, in which case the continuity of $L_1$ implies the continuity of $L_2$, or $L_1(\nu)=L_2(\nu)$, where the situation is trivial. If $L_3(\mu)>L_4(\mu)$, considering $\wedge_3A$, through a similar analysis, we can not only obtain the continuity of $L_1(\nu) + L_2(\nu) + L_3(\nu)$, but also further establish the continuity of $ L_1,L_2,\text{ and }L_3$ individually. Moreover, for $L_k(\mu)>L_{k+1}(\mu)$ where $1\leq k\leq d$, the continuity of $L_k$ can be obtained by repeating the above process.
 
\begin{remark}
	On the other hand, if we assume that the fiber action $A$ is invertible. Then by applying properties of the Wasserstein distance, we can conclude that $$P(X)\ni\nu\mapsto L_1(\nu)+\cdots+L_d(\nu)=\int_X\log|\det(A(x))|d\nu(x)$$ is continuous everywhere since continuous functions can be approximated by Lipschitz continuous functions. It naturally enables us to establish the continuity of all Lyapunov exponents regardless of whether a gap exists or not. The assumption for invertibility seems crucial. As is well known, a set of measure zero can become a set of positive measure under slight perturbations of the measure. It could cause the Lyapunov exponents to drop directly from a finite number to negative infinity. 
\end{remark}
\begin{remark}
	In fact, for $1\leq k\leq d$, by applying Lemma \ref{lem1} in each step of the iteration, we have $$L_1^{(n_k)}(\nu)-L_2^{(n_k)}(\nu)>(\kappa(\mu)-2\theta_k-3\varepsilon)(1-r_3(n_k)).$$ Then taking the limit as $k\to\infty$, it is evident that
	$$L_1(\nu)-L_2(\nu)\geq \kappa(\mu)-46\varepsilon-3\varepsilon>L_1(\mu)-L_2(\mu)-50\varepsilon,$$    
	which proves the lower semicontinuity of $L_1(\nu)-L_2(\nu)$.
\end{remark}

From the estimate \eqref{g10},  for an increasing sequence $n=n_k,k>0$  we have $$ |L_1^{(n)}(\nu)-L_1^{(n)}(\mu)|<\varepsilon.$$ 
\par Now we consider the issue on intervals of scales rather than individual scales. This further illustrates that the conclusion is valid at all large enough scales $n$. 
\par Denote $\mathcal{N}_0=[m_1,e^{m_1}]:=[n_0^-,n_0^+]$. To begin with, the base step of the inductive procedure holds for not only a single scale $n_0$  but also $\mathcal{N}_0$ with $m_1$ from Proposition \ref{prop2}. Let $\Gamma(n)=n^{1+}$, then we can define inductively
\begin{align*}
    & \mathcal{N}_1=[\Gamma(n_0^-),\Gamma(n_0^+)]:=[n_1^-,n_1^+],\\
    &\mathcal{N}_{k+1}=[\Gamma(n_k^-),\Gamma(n_k^+)]:=[n_{k+1}^-,n_{k+1}^+].
\end{align*}
If $n=n_{k+1}\in\mathcal{N}$, then there is $n=n_{k+1}\asymp\Gamma(n_k)=n_k^{1+}$ such that $$|L_1^{(n_{k+1})}(\nu)-L_1^{(n_{k+1})}(\mu)|<\theta_{k+1}.$$
Note that $\Gamma(n)$ is an increasing function and $\mathcal{N}_k$ and $\mathcal{N}_{k+1}$ overlap, so we can continue this procedure. 
\par Moreover,we can apply Lemma \ref{lem1} for all large enough scales $n$. Therefore, we deduce that the following uniform, finite scale statement holds. Note that this is more general than the result of Proposition \ref{prop2}.
\begin{lemma}
    Given a cocycle $(\mu,A)\in\mathscr{C}$ with $L_1(\mu)>L_2(\mu)$ and $0<\varepsilon<\kappa(\mu)/50$, there are $\delta=\delta(\mu,\varepsilon)>0, m=m(\mu,\varepsilon)\in\mathbb{N}$ and $c=c(\mu,\varepsilon)>0$ such that for all $n\geq m$ and for all $\nu\in P(X)$ with $W_1(\mu,\nu)<\delta$ we have:
    \begin{align}
        & |L_1^{(n)}(\nu)-L_1^{(n)}(\mu)|<\varepsilon,\label{in1}\\
        & \frac{1}{n}\log gr(A^{(n)}(x))>\kappa(\mu)-5\varepsilon.\label{in2}
    \end{align}
    for all $x$ outside a set of $\nu$-measure $<r(n)$.
\end{lemma}

\subsection{Modulus of Continuity}\label{sec:mod}
For the purpose of getting H\"older continuity of Lyapunov blocks w.r.t. measure, we will specify the following proposition, it gives the rate of convergence of the finite scale exponents $L_1^{(n)}(\nu)$ to the first Lyapunov exponent $L_1(\nu)$ and it gives an estimate on the proximity of these finite scale Lyapunov exponents at different scales. Let $\phi$ be the inverse of $\psi$ defined in Proposition \ref{prop3} and for every integer $n\in\mathbb{N}$, denote $n++:=\lfloor\psi(n)\rfloor=\lfloor n[r(n)]^{-(\frac{1}{2}-)}\rfloor,n--:=\lfloor\phi(n)\rfloor$, so $(n++)--\asymp n$.
\begin{proposition}[uniform speed of convergence]\label{prop4}
     Let $(\mu,A)\in\mathscr{C}$ be a cocycle for which $L_1(\mu)>L_2(\mu)$. There are $\delta=\delta(\mu)>0,C_2=C_2(\mu),m_4=m_4(\mu)\in\mathbb{N},c_3=c_3(\mu)>0$ such that for all $n\geq m_4$ and for all $\nu\in P(X)$ with $W_1(\nu,\mu)<\delta$, we have
     \begin{align}
        L_1^{(n)}(\nu)-L_1(\nu)<C_2\frac{\phi(n)}{n}\leq [r_3(n--)]^{\frac{1}{2}-}, \label{u1}\\
        |L_1^{(n++)}(\nu)+L_1^{(n)}(\nu)-2L_1^{(2n)}(\nu)|<C_2\frac{n}{n++}\leq [r_3(n)]^{\frac{1}{2}-}.\label{u2}
    \end{align}
\end{proposition}
\begin{proof}
    Fix $0<\varepsilon<\kappa(\mu)/50$. We consider parameters $\delta\leq\min\{\delta_1,\delta_2\},m\geq\max\{m_1,m_2\} $, and $c_3$ is the notation established in the proof of Theorem \ref{thm2}. Additionally, $C_1,C_2$ are the constants in Proposition \ref{prop2} and Proposition \ref{prop3} respectively.
    \par Pick $n_0^-\in\mathbb{N}$ large enough such that Proposition \ref{prop2} and Proposition \ref{prop3} hold for all $n\geq n_0^-$, and $L_1^{(n)}(\mu)-L_1^{(2n)}(\mu)<\varepsilon$ for $n\geq n_0^-$. Moreover, assume also $n_0^-$ to be large enough so that $e^{-C_12e^{n_0^-}}<\delta,r_3(n_0^-)<\varepsilon,n^{0+}\ll [r_3(n)]^{-(\frac{1}{2}-)}$ for $n\geq n_0^-$, $C_2\frac{n_0^-}{(n_0^-)^{1+}}<\varepsilon$ and $n_0^-[r_3(n_0^-)]^{\frac{1}{2}-}<e^{n_0^-}$.
    \par Next, we will discuss the problem on intervals of scales, rather than on individual scales.
    \par Set $n_0^+:=\lfloor e^{n_0^-}\rfloor,\mathcal{N}_0:=[n_0^-,n_0^+]$ and $\delta:=e^{-C_12n_0^+}$. Then for all measures $\nu$ with $W_1(\nu,\mu)<\delta$, and for all $n_0\in\mathcal{N}_0$, since $\delta=e^{-2C_1n_0^+}\leq e^{-2C_1n_0}<e^{-C_1n_0}$, we can apply Proposition \ref{prop2} at scales $2n_0$ and $n_0$ so that the assumptions in the inductive step Proposition \ref{prop3} hold for every $n_0\in\mathcal{N}_0$.
    \par Define inductively the intervals of scales for all $k\geq 0$:
    \begin{align*}
        &\mathcal{N}_1=[\psi(n_0^-),\psi(n_0^+)]=[n_1^-,n_1^+],\\
        &\mathcal{N}_{k+1}=[\psi(n_k^-),\psi(n_k^+)]=[n_{k+1}^-,n_{k+1}^+].
    \end{align*}
    \par Therefore, if $n_1\in\mathcal{N}_1$, then for some $n_0\in\mathcal{N}_0$,
    $$n_1=n_0\lfloor [r_3(n_0)]^{-(\frac{1}{2}-)}\rfloor\asymp n_0[r_3(n_0)]^{-(\frac{1}{2}-)}=\psi(n_0).$$
    \par Apply Proposition \ref{prop3} and obtain
    \begin{align*}
        &L_1^{(n_1)}(\nu)-L_1^{(2n_1)}(\nu)<C_2\frac{n_0}{n_1}\asymp C_2\frac{\phi(n_1)}{n_1},\\
        &|L_1^{(n_1)}(\nu)+L_1{(n_0)}(\nu)-L_1^{(2n_0)}(\nu)|<C_2\frac{n_0}{n_1}\asymp C_2[r_3(n_0)]^{\frac{1}{2}-},
    \end{align*}
    since $n_1\asymp\psi(n_0)$, so $\phi(n_1)\asymp n_0$.
    \par Note that $\psi(n)$ is an increasing function and $\mathcal{N}_k$ and $\mathcal{N}_{k+1}$ overlap, so continue inductively, if $n_{k+1}\in\mathcal{N}_{k+1}$, then there is $n_k\in\mathcal{N}_k$ such that 
    $$n_k++=n_{k+1}=n_k\lfloor [r_3(n_k)]^{-(\frac{1}{2}-)}\rfloor\asymp n_0[r_3(n_k)]^{-(\frac{1}{2}-)}=\psi(n_k).$$
    Then 
    \begin{align}
       &L_1^{(n_{k+1})}(\nu)-L_1^{(2n_{k+1})}(\nu)<C_2\frac{n_k}{n_{k+1}}\asymp C_2\frac{\phi(n_{k+1})}{n_{k+1}}\label{u3}\\
       &|L_1^{(n_k++)}(\nu)+L_1^{(n_k)}(\nu)-2L_1^{(2n_k)}(\nu)|<C_2\frac{n_k}{n_k++}\asymp C_2[r_3(n_k)]^{\frac{1}{2}-}.\label{u4}
    \end{align}
    \par In fact, arbitrary adjacent intervals $\mathcal{N}_k\text{ and }\mathcal{N}_{k+1}$ for $k\geq 0$ are intersecting since $\psi$ is increasing.   Hence, for all $n\geq n_1^-$,
    \begin{align}
       &L_1^{(n)}(\nu)-L_1^{(2n)}(\nu)<C_2\frac{\phi(n)}{n}\leq [r_3(n--)]^{\frac{1}{2}-},\label{u5}\\
       &|L_1^{(n++)}(\nu)+L_1^{(n)}(\nu)-2L_1^{(2n)}(\nu)|<C_2\frac{n}{n++}\leq [r_3(n)]^{\frac{1}{2}-}.\label{u6}
    \end{align}
    are proved.
    \par Note that \eqref{u5} implies that for all $k\geq 0$,
    \begin{equation}\label{u7}
        L_1^{(2^kn)}(\nu)-L_1^{(2^{k+1}n)}(\nu)<C_2\frac{\phi(2^kn)}{2^kn}.
    \end{equation}
    \par Since for all deviation measure functions $r(n)\in\mathcal{J}$, the corresponding map $\phi=\phi_r$ satisfies
    $$\limsup_{n\to\infty}\frac{\phi(2n)}{\phi(n)}\leq 2\beta<2,$$
    where $0<\beta<1$. Then $\phi(2n)\leq 2\beta\phi(n)$ for large enough $n$.
    \par Therefore, $\frac{\phi(2^kn)}{2^kn}\leq \beta^k\frac{\phi(n)}{n}$ holds for all $k\geq 0$.
    \par Sum up \eqref{u7} from $k=0$ to $\infty$,
    \begin{align*}
    	L_1^{(n)}(\nu)-L_1(\nu)&=\sum_{k=0}^{\infty}\big(L_1^{(2^kn)}(\nu)-L_1^{(2^{k+1}n)}(\nu)\big)\\
    	&<\sum_{k=0}^{\infty}\frac{\phi(2^kn)}{2^kn}\leq \frac{\phi(n)}{n}\cdot\frac{1}{1-\beta}\leq C_2\frac{\phi(n)}{n}.
    \end{align*}
    which implies that \eqref{u1} is proven.
\end{proof}

\begin{theorem}[modulus of continuity]\label{thm3}
    Let $A\in\mathscr{C}$ be a measurable cocycle for which $L_1(\mu)>L_2(\mu)$. There are $\delta=\delta(\mu)>0,c_4=c_4(\mu)>0,a=a(\mu)>0$ such that if we define the modulus of continuity function $\omega(h):=[r_4(a\log\frac{1}{h})]^{\frac{1}{2}-}$, then for all measures $\nu_i\in P(X)(i=1,2)$ with $W_1(\nu_i,\mu)<\delta$, we have
     \begin{equation}\label{m1}
         |L_1(\nu_1)-L_1(\nu_2)|\leq\omega(W_1(\nu_1,\nu_2)).
     \end{equation}
\end{theorem}
\begin{proof}
    Choose $\delta$ be less than the size of the neighborhood of $\mu$ from Proposition \ref{prop4}, $m\geq\max\{m_1,m_4\}$ that $m_1,m_4$ from Proposition \ref{prop2} and \ref{prop4}, $r_4(n)=e^{-c_4n}\geq r_1(n)+r_3(n)$ be the corresponding deviation measure function such that the Proposition \ref{prop2} and \ref{prop4} apply for all $\nu_i$ with $W_1(\nu_i,\mu)<\delta$ and for all $n\geq\bar{n}_0$.
    \par Let $C_1=C_1(\mu)$ be the constant from Proposition \ref{prop2}.
    \par Let $W_1(\nu_1,\nu_2):=h$ and set $n:=\lfloor\frac{1}{2C_1}\log\frac{1}{h}\rfloor\in\mathbb{N}$, then $e^{-C_14n}\leq h\leq e^{-C_12n}$, so $W_1(\nu_1,\nu_2)=h<e^{-C_12n}$ and $n\geq m$.
    \par Applying Proposition \ref{prop2} to $\nu_1,\nu_2$ at scales $2n$ and $n$, we get
    \begin{align}
         &|L_1^{(n)}(\nu_1)-L_1^{(n)}(\nu_2)|<r_4(n)<[r_4(n)]^{\frac{1}{2}-}\label{m3}\\
         &|L_1^{(2n)}(\nu_1)-L_1^{(2n)}(\nu_2)|<r_4(n)<[r_4(n)]^{\frac{1}{2}-}\label{m4}
    \end{align}
    \par Applying Proposition \ref{prop4} to $\nu_1,\nu_2$ at scales $n$, we have
     \begin{align}
        L_1^{(n++)}(\nu_i)-L_1(\nu_i)<[r_4(n)]^{\frac{1}{2}-}\label{m5}\\
        |L_1^{(n++)}(\nu_i)+L_1^{(n)}(\nu_i)-2L_1^{(2n)}(\nu_i)|<[r_4(n)]^{\frac{1}{2}-}.\label{m6}
    \end{align}
    \par Combine inequalities \eqref{m3}-\eqref{m6}, and we can obtain
    \begin{equation*}
       \begin{split}
           |L_1(\nu_1)-L_1(\nu_2)|=&\Biggl|L_1^{(n++)}(\nu_2)-L_1(\nu_2)-\big(L_1^{(n++)}(\nu_1)-L_1(\nu_1)\big)\\&+\big(L_1^{(n++)}(\nu_1)+L_1^{(n)}(\nu_1)-2L_1^{(2n)}(\nu_1)\big)\\&-\big(L_1^{(n++)}(\nu_2)+L_1^{(n)}(\nu_2)-2L_1^{(2n)}(\nu_2)\big)\\&+\big(L_1^{(n)}(\nu_2)-L_1^{(n)}(\nu_1)\big)+2\big(L_1^{(2n)}(\nu_1)-L_1^{(2n)}(\nu_2)\big)\Biggr|\\&\leq 7[r_4(n)]^{\frac{1}{2}-}\leq h^{-\frac{c_4}{C_1}\cdot(\frac{1}{4}-)}=\omega(h)=\omega(W_1(\nu_1,\nu_2)).
       \end{split}
    \end{equation*}
\end{proof}
\begin{remark}
	More generally, if for some $1\leq k\leq d$, the cocycle $(\mu,A)$ has the Lyapunov spectrum gap $L_k(\mu)>L_{k+1}(\mu)$, then taking exterior powers, we can conclude that
	\begin{align*}
		L_1(\wedge_kA)&=(L_1+\cdots +L_{k-1}+L_k)(\mu)\\&>(L_1+\cdots +L_{k-1}+L_{k+1})(\mu)=L_2(\wedge_kA).
	\end{align*}
	Then the map $\Lambda_k:=L_1+\cdots+L_k$ satisfies
	$$|\Lambda_k(\nu_1)-\Lambda_k(\nu_2)|=|L_1(\wedge_k\nu_1)-L_1(\wedge_k\nu_2)|\leq\omega(W_1(\nu_1,\nu_2)).$$
\end{remark}
\begin{remark}
    If the deviation set measure functions $r(n)$ decay at most sub-exponentially fast, we can only achieve weak H\"older continuity of the Lyapunov exponents, which is significantly weaker than our results.
\end{remark}

\subsection{Further Discussion on the Joint Continuity Theorem}
Having established the ACT in measure of the Lyapunov exponents in the previous sections, we now turn to extend our result in the context of perturbing both measure and fiber.
\par Consider that the space of matrix-valued functions is endowed with the $C_0$ norm while the probability measure space is endowed with the $W_1$ distance. When the measure and fiber versions of the uniform base and fiber-LDT estimates hold simultaneously, combining our main Theorem \ref{thm1} and Theorem 3.1 of Duarte-Klein \cite{duarte2016lyapunov}, we can directly derive the H\"older continuity of the first Lyapunov exponent $L_1$ with respect to the measure and the fiber, respectively. Furthermore, for any $(\nu, B)$ sufficiently close to $(\mu, A)$, by applying the triangle inequality, we can decompose as follows:  
$$|L_1(\nu,B)-L_1(\mu,A)| \leq |L_1(\nu,B)-L_1(\nu,A)| + |L_1(\nu,A)-L_1(\mu,A)|,$$ 
which naturally leads to the joint Hölder continuity of $L_1$ as a function of the measure and the fiber as follows.
\begin{theorem}
	Given a continuous cocycle $(\mu,A)\in\mathscr{C}$ such that\\
	(1) every $\xi\in \mathcal{C}^\alpha(X)(0<\alpha\leq 1)$ satisfies a uniform base-LDT estimate.\\
	(2) every cocycle $(\mu,A)\in\mathscr{C}$ for which $L_1(\mu,A)>L_2(\mu,A)$ satisfies a uniform fiber-LDT estimate (in measure and fiber).\par
	There exists a neighborhood $\mathcal{V}$ of $(\mu,A)$ such that for all $(\nu,B)\in\mathcal{V}$, the first Lyapunov exponent $L_1(\nu,B)$ is H\"older continuous on the measure and fiber. 
\end{theorem}

\section{Application to the product of random matrices}\label{apply}
Let $\Sigma=\mathrm{SL}(d,\mathbb{R})$ be a Polish space and $(\Sigma,\mathcal{B},\mu)$ be a probability space. Consider the product space $X=\Sigma^{\mathbb{Z}}$, endowed with the product $\sigma$-algebra $\mathcal{F}=\mathcal{B}^{\mathbb{Z}}$ and the product
measure $\mu^{\mathbb{Z}}$ where $\mu\in \mathrm{Prob} (\mathrm{SL}(d,\mathbb{R}))$. The two-sided shift is $\sigma:X\to X$ defined by $$\sigma((x_n)_n)=(x_{n+1})_n,$$ which is a homeomorphism and preserves the measure $\mu^{\mathbb{Z}}$. We call the MPDS $(X,\mathcal{F},\mu,\sigma)$ a Bernoulli shift. Note that every Bernoulli shift is always strongly mixing.
\par Consider a continuous function $A:X\to\mathrm{SL}(d,\mathbb{R})$ and $\mathrm{supp}(\mu)$ being contained in a compact set. The measure $\mu$ and matrix-valued function $A$ determine a cocycle $(\mu,A)$, and iterates of the cocycle are given by $$A^{(n)}(x)=A(\sigma^{n-1}x)\cdots A(\sigma x)A(x),$$ where $x=(x_n)_{n\in\mathbb{Z}}\in X$.
\par Let $\{A_n\}_{n\geq 1}$ be a sequence of i.i.d. $d\times d$ real invertible random matrices with determinant of 1 on the cocycle $(\mu,A)$. Denote the product by $P_n=A_1A_2\cdots A_n$. As stated in Theorem 2 of Furstenberg and Kesten \cite{furstenberg1960products}, if $\mathbb{E}(\log^+\|A_1\|)<\infty$, then for $\mu$-a.e. $x\in X$, $$\lim_{n\to\infty}\frac{1}{n}\log\|P_n(x)\|=L_1(\mu)$$ with $L_1(\mu)$ being called the first Lyapunov exponent of the sequence $\{A_n\}_{n\geq 1}$.
\par Within our framework, it is possible to derive the following uniform (in measure) base and fiber LDT estimates and establish the continuity of Lyapunov exponents for random matrix products. Duarte-Klein \cite{duarte2016lyapunov} have demonstrated the uniform LDT estimates in fiber for cocycles over mixing Markov chains, which represents a more challenging case. Drawing on their model and methodology, we can easily derive the results of measure version as follows.
\begin{theorem}\label{base}\cite{duarte2016lyapunov}
	For an observable $\xi\in \mathcal{C}^\alpha(X)(0<\alpha\leq 1)$ and all $\varepsilon>0$, there exists a neighborhood $\mathcal{V}_\mu$ of $\mu$ and there exist $C=C(\xi,\varepsilon,\mu)<\infty,c=c(\xi,\varepsilon,\mu)>0$ such that for all $\nu\in\mathcal{V}_\mu$ and $n\in\mathbb{N}$,
	$$\nu^{\mathbb{Z}}\left\{x\in X:\left|\frac{1}{n}\sum_{j=0}^{n-1}\xi\circ \sigma^j(x)-\mathbb{E}_\mu(\xi)\right|>\varepsilon\right\}\leq r(n).$$
\end{theorem}
\begin{theorem}\label{fiber}\cite{duarte2016lyapunov}
	Assume that $\mu^{\mathbb{Z}}$ is quasi-irreducible and $L_1(\mu)>L_2(\mu)$. Then for all $\varepsilon>0$, there exist a neighborhood $\mathcal{V}_\mu$ of $\mu,C=C(A,\mu,\varepsilon)<\infty \text{ and }c=c(A,\mu,\varepsilon)>0$ such that for all $\nu\in\mathcal{V}_\mu$ and $n\in\mathbb{N}$, 
	$$\nu^{\mathbb{Z}}\left\{x\in X:\left|\frac{1}{n}\log\|A^{(n)}(x)\|-L_1(\nu)\right|>\varepsilon\right\}\leq r(n).$$
\end{theorem}
\begin{remark}
	To maintain consistency with the previous work, we still use \(r(n)\) to control the measure of the deviation set. In fact, as referenced in Duarte-Klein \cite{duarte2016lyapunov}, we can take \(r(n)=O(\varepsilon^2)\).
\end{remark}
\begin{theorem}\label{application}
    Assume that the measure $\mu$ is quasi-irreducible and $L_1(\mu)>L_2(\mu)$, then the first Lyapunov exponent of random matrix products is H\"older continuous with respect to measure in Wasserstein distance. 
\end{theorem}
\begin{proof}
	This is a straightforward application of the theorem established earlier. It is well-known that the Bernoulli shift is always strongly mixing. Referring to Section 5.3.1 in \cite{duarte2016lyapunov}, we can prove that the uniform base LDT estimates namely Theorem \ref{base} holds with respect to the measure. Moreover, we use the assumptions of quasi-irreducibility and \(L_1(\mu) > L_2(\mu)\), proving uniform fiber LDT estimates Theorem \ref{fiber} for the measure. The detailed arguments are the same as that in Section 5.3.2 in \cite{duarte2016lyapunov}. Consequently, by directly applying the measure version of the ACT established in Wasserstein distance, we obtain the H\"older continuity of the first Lyapunov exponent for product of random matrices.
\end{proof}
\begin{remark}
	Since the total variation distance is always greater than the Wasserstein distance, we cannot directly obtain H\"older continuity from Lipschitz continuity proved in Baraviera-Duarte \cite{baraviera2019approximating}.
\end{remark}
\par In particular, given an i.i.d. sequence of random variables $\omega=\{\omega_n\}_{n\in\mathbb{Z}}$ which represents the random potential, consider the one-dimensional Schr\"odinger operator $H$ defined by
$$(H\psi)n=-\psi_{n+1}-\psi_{n-1}+\omega_n\psi_n,\quad \forall n\in\mathbb{Z},$$
for all $\psi=\{\psi_n\}_{n\in\mathbb{Z}}\in l^2(\mathbb{Z})$. Then for the energy parameter $E$, the corresponding transfer matrix $S_E(\omega)$ is given by 
$$S_E(\omega)=\begin{pmatrix}
	\omega-E&-1\\1&0 
\end{pmatrix}.$$ 
\par Given $\rho\in\mathrm{Prob}_c(\mathbb{R})$, let 
$$\mu_E:=\int_{\mathbb{R}}\delta_{S_E(\omega)}d\rho(\omega)\in\mathrm{Prob}(\mathrm{SL}(2,\mathbb{R})),$$
Hence, the random Schr\"odinger cocycle is driven by the MPDS $(\mathrm{SL}(2,\mathbb{R})^{\mathbb{Z}},\sigma,\mu_E^{\mathbb{Z}},\mathcal{F})$, and the iterates of the cocycle are defined by $$S_E^{(n)}(\omega)=S_E(\omega_{n-1})\cdots S_E(\omega_1)S_E(\omega_0).$$
We denote the first Lyapunov exponent of the cocycle by
$$L_1(\mu_E)=\lim_{n\to\infty}\frac{1}{n}\log\|S_E^{(n)}(\omega)\|, \text{ for } \mu_E \text{-a.e. } \omega.$$
\par For the proof of Corollary \ref{cor}, by Furstenberg \cite{furstenberg1963noncommuting}, it follows that when measure $
\rho$ is not a single Dirac, strong irreducibility and $L_1(\mu_E)>L_2(\mu_E)$ hold. So for the measure $\nu$ in a neighborhood of $\mu_E$, we can directly apply Theorem \ref{application} to conclude that the first Lyapunov exponent is H\"older continuous with respect to measure $\nu$ in Wasserstein distance, which implies the continuity on energy $E$.
\par In fact, in the field of mathematical physics, people are particularly concerned with the continuity of Lyapunov exponents with respect to energy $E$. The continuity of Lyapunov exponents is not only of dynamical importance but also closely tied to the properties of the Integrated Density of States (IDS) in the spectral theory of Schr\"odinger operators. Specifically, (weak) H\"older continuity of the Lyapunov exponents with respect to energy $E$ is equivalent to the same regularity of the IDS by the famous Thouless formula. This equivalence offers profound insights into the dynamical stability and spectral properties of the system, and further establishes a robust mathematical framework for analyzing complex physical models.
\section*{Acknowledgments}
The authors would like to thank Pedro Duarte and Silvius Klein for their useful suggestions to improve various aspects of this manuscript. This work was supported by NSFC grant 12301231.
%\bibliographystyle{amsplain}
%\bibliography{reference}

\begin{thebibliography}{99}
	
	\bibitem{arnold1995random}
	L.~Arnold.
	\newblock {\em {R}andom dynamical systems}.
	\newblock Springer Berlin, Heidelberg, 1998.
	
	\bibitem{avila2023continuity}
	A.~Avila, A.~Eskin, and M.~Viana.
	\newblock {C}ontinuity of the lyapunov exponents of random matrix products.
	\newblock {\em arXiv preprint arXiv:2305.06009}, 2023.
	
	\bibitem{avila2010extremal}
	A.~Avila and M.~Viana.
	\newblock {E}xtremal lyapunov exponents: an invariance principle and
	applications.
	\newblock {\em Invent. Math.}, 181(1):115--178, 2010.
	
	\bibitem{avila2014complex}
	A.~Avila, S.~Jitomirskaya, and C.~Sadel.
	\newblock {C}omplex one-frequency cocycles.
	\newblock {\em J. Eur. Math. Soc. (JEMS)}, 16(9):1915--1935, 2014.
	
	\bibitem{backes2018continuity}
	L.~Backes, A.~Brown, and C.~Butler.
	\newblock {C}ontinuity of lyapunov exponents for cocycles with invariant
	holonomies.
	\newblock {\em J. Mod. Dyn.}, 12:223--260, 2018.
	
	\bibitem{baraviera2019approximating}
	A.~Baraviera and P.~Duarte.
	\newblock {A}pproximating lyapunov exponents and stationary measures.
	\newblock {\em J. Dynam. Differential Equations}, 31:25--48, 2019.
	
	\bibitem{bezerra2019random}
	J.~Bezerra and M.~Poletti.
	\newblock {R}andom product of quasi-periodic cocycles.
	\newblock {\em arXiv preprint arXiv:1907.00815}, 2019.
	
	\bibitem{BOCHI2002}
	J.~Bochi.
	\newblock {G}enericity of zero lyapunov exponents.
	\newblock {\em Ergodic Theory Dynam. Systems}, 22(6):1667--1696, 2002.
	
	\bibitem{bocker2017continuity}
	C.~Bocker-Neto and M.~Viana.
	\newblock {C}ontinuity of lyapunov exponents for random two-dimensional
	matrices.
	\newblock {\em Ergodic Theory Dynam. Systems}, 37(5):1413--1442, 2017.
	
	\bibitem{bourgain2002continuity}
	J.~Bourgain and S.~Jitomirskaya.
	\newblock {C}ontinuity of the lyapunov exponent for quasiperiodic operators
	with analytic potential.
	\newblock {\em J. Stat. Phys.}, 108:1203--1218, 2002.
	
	\bibitem{cai2022mixed}
	A.~Cai, P.~Duarte, and S.~Klein.
	\newblock {M}ixed random-quasiperiodic cocycles.
	\newblock {\em Bull. Braz. Math. Soc. (N.S.)}, 53(4):1469--1497, 2022.
	
	\bibitem{cai2023furstenberg}
	A.~Cai, P.~Duarte, and S.~Klein.
	\newblock {F}urstenberg theory of mixed random-quasiperiodic cocycles.
	\newblock {\em Comm. Math. Phys.}, 402(1):447--487, 2023.
	
	\bibitem{cramer1994collected}
	H.~Cram{\'e}r.
	\newblock {\em {C}ollected Works II''}.
	\newblock Springer Berlin, Heidelberg, 1994.
	
	\bibitem{dembo2009large}
	A.~Dembo.
	\newblock {\em {L}arge deviations techniques and applications''}.
	\newblock Springer Berlin, Heidelberg, 2009.
	
	\bibitem{donsker1975asymptotic1}
	M.~D. Donsker and S.~R.~S. Varadhan.
	\newblock {A}symptotic evaluation of certain markov process expectations for
	large time, i.
	\newblock {\em Comm. Pure Appl. Math.}, 28(1):1--47, 1975.
	
	\bibitem{donsker1975asymptotic2}
	M.~D. Donsker and S.~R.~S. Varadhan.
	\newblock {A}symptotic evaluation of certain markov process expectations for
	large time, ii.
	\newblock {\em Comm. Pure Appl. Math.}, 28(2):279--301, 1975.
	
	\bibitem{donsker1976asymptotic3}
	M.~D. Donsker and S.~R.~S. Varadhan.
	\newblock {A}symptotic evaluation of certain markov process expectations for
	large time, iii.
	\newblock {\em Comm. Pure Appl. Math.}, 29(4):389--461, 1976.
	
	\bibitem{donsker1983asymptotic4}
	M.~D. Donsker and S.~R.~S. Varadhan.
	\newblock {A}symptotic evaluation of certain markov process expectations for
	large time, iv.
	\newblock {\em Comm. Pure Appl. Math.}, 36(2):183--212, 1983.
	
	\bibitem{duartefreijo2024continuity}
	P.~Duarte and C.~Freijo.
	\newblock {C}ontinuity of the lyapunov exponents of non-invertible random
	cocycles with constant rank.
	\newblock {\em Nonlinearity}, 37(10):105008, 2024.
	
	\bibitem{duarte2014continuity}
	P.~Duarte and S.~Klein.
	\newblock {C}ontinuity of the lyapunov exponents for quasiperiodic
	cocycles.
	\newblock {\em Comm. Math. Phys.}, 332:1113--1166, 2014.
	
	\bibitem{duarte2016lyapunov}
	P.~Duarte and S.~Klein.
	\newblock {L}yapunov exponents of linear cocycles: Continuity via large
	deviations.
	\newblock {\em Atlantis Studies in Dynamical Systems}, 3, 2016.
	
	\bibitem{duarte2020large}
	P.~Duarte and S.~Klein.
	\newblock {L}arge deviations for products of random two dimensional
	matrices.
	\newblock {\em Comm. Math. Phys.}, 375:2191--2257, 2020.
	
	\bibitem{hislop2020dependence}
	P.~D.~ Hislop and C.~A.~ Marx. \emph{Dependence of the density of states
		on the probability distribution for discrete random schr{\"o}dinger
		operators}, Int. Math. Res. Not., 17: 5279--5341, 2020.
	
	\bibitem{pinto2021holder}
	M.~Durães.
	\newblock {\em {H}\"older continuity for Lyapunov exponents of random linear
		cocycles}.
	\newblock PhD thesis, PUC-Rio, 2021.
	\newblock PhD dissertation.
	
	\bibitem{furstenberg1963noncommuting}
	H.~Furstenberg.
	\newblock {N}oncommuting random products.
	\newblock {\em Trans. Amer. Math. Soc.}, 108(3):377--428, 1963.
	
	\bibitem{furstenberg1960products}
	H.~Furstenberg and H.~Kesten.
	\newblock {P}roducts of random matrices.
	\newblock {\em Ann. Math. Stat.}, 31(2):457--469, 1960.
	
	\bibitem{furstenbergkifer1983}
	H.~Furstenberg and Y.~Kifer.
	\newblock {R}andom matrix products and measures on projective spaces.
	\newblock {\em Israel J. Math.}, 46:12--32, 1983.
	
	\bibitem{goldstein2001holder}
	M.~Goldstein and W.~Schlag.
	\newblock {H}{\"o}lder continuity of the integrated density of states for
	quasi-periodic schr{\"o}dinger equations and averages of shifts of
	subharmonic functions.
	\newblock {\em Ann. of Math. (2)}, 154:155--203, 2001.
	
	\bibitem{jitomirskaya2009continuity}
	S.~Jitomirskaya, D.~A. Koslover, and M.~S. Schulteis.
	\newblock {C}ontinuity of the lyapunov exponent for analytic quasiperiodic
	cocycles.
	\newblock {\em Ergodic Theory Dynam. Systems}, 29(6):1881--1905, 2009.
	
	\bibitem{kifer2012ergodic}
	Y.~Kifer.
	\newblock {\em {E}rgodic theory of random transformations}, volume~10.
	\newblock Springer Science \& Business Media, 2012.
	
	\bibitem{le1989regularite}
	\'E. Le~Page.
	\newblock {R}\'egularit\'e du plus grand exposant caract\'eristique des
	produits de matrices al\'eatoires ind\'ependantes et applications.
	\newblock {\em Ann. Inst. Henri Poincar\'e Probab. Stat.}, 25(2):109--142,
	1989.
	
	\bibitem{malheiro2015lyapunov}
	E.~C. Malheiro and M.~Viana.
	\newblock {L}yapunov exponents of linear cocycles over markov shifts.
	\newblock {\em Stoch. Dyn.}, 15(03):1550020, 2015.
	
	\bibitem{rassoul2015course}
	F.~Rassoul-Agha and T.~Sepp{\"a}l{\"a}inen.
	\newblock {\em {A} course on large deviations with an introduction to Gibbs
		measures}, volume 162.
	\newblock American Mathematical Society, 2015.
	
	\bibitem{sanchez2018lyapunov}
	A.~S{\'a}nchez and M.~Viana.
	\newblock {L}yapunov exponents of probability distributions with non-compact
	support.
	\newblock {\em arXiv preprint arXiv:1810.03061}, 2018.
	
	\bibitem{schlag2013regularity}
	W.~Schlag.
	\newblock {R}egularity and convergence rates for the lyapunov exponents of
	linear cocycles.
	\newblock {\em J. Mod. Dyn.}, 7(4):619--637, 2013.
	
	\bibitem{shamis2021continuity}
	M.~Shamis. \emph{On the continuity of the integrated density of states in the
		disorder}, Int. Math. Res. Not., 22: 17304--17315, 2021.
	
	\bibitem{tall2020moduli}
	E.~Y.~Tall and M.~Viana. \emph{Moduli of continuity for the lyapunov
		exponents of random GL(2)-cocycles}, Trans. Amer. Math. Soc., 373(2): 1343--1383, 2020.
	
	\bibitem{varadhan1966asymptotic}
	S.~R.~S. Varadhan.
	\newblock {A}symptotic probabilities and differential equations.
	\newblock {\em Comm. Pure Appl. Math.}, 19(3):261--286, 1966.
	
	\bibitem{viana2016foundations}
	M.~Viana.
	\newblock {\em {F}oundations of Ergodic Theory}, volume 151.
	\newblock Cambridge University Press, 2016.
	
	\bibitem{villani2009optimal}
	C.~Villani.
	\newblock {\em {O}ptimal transport: old and new}, volume 338.
	\newblock Springer Berlin, Heidelberg, 2009.
	
\end{thebibliography}
%\nocite{*}

\end{document}